%% file: main.tex
\documentclass[a4paper,11pt,reqno]{amsart}
\pdfoutput=1

\usepackage[english]{babel}
\usepackage{xcolor}
\usepackage[pdfencoding=unicode, hidelinks]{hyperref}
\hypersetup{
	colorlinks=false,
	linkcolor={red},
	citecolor={blue},
}
\usepackage{crossreftools}
\usepackage{etoolbox}
\cslet{blx@noerroretextools}\empty 
\usepackage{mathtools}
\usepackage{mathrsfs}
\usepackage{amsthm}
\usepackage[nameinlink]{cleveref}
\usepackage{amsfonts}
\usepackage{amssymb}
\usepackage{amsmath}
\usepackage{autonum}
\usepackage[maxnames=50, sorting=ayt, doi=false, isbn=false, url=false, giveninits=true, date=year]{biblatex}
\AtEveryBibitem{\clearfield{number}}
\renewbibmacro{in:}{}
\usepackage{tikz}
\usepackage{orcidlink}
\usepackage{float}
\usepackage{stmaryrd}
\usepackage{graphics}
\usepackage{graphicx}
\usepackage{subfig}
\usepackage{cases}

\usepackage[normalem]{ulem}

\DeclareSortingTemplate{ayt}{
  \sort{
    \field{presort}
  }
  \sort[final]{
    \field{sortkey}
  }
  \sort{
    \field{sortname}
    \field{author}
  }
  \sort{
    \field{sortyear}
    \field{year}
    \literal{9999}
  }
  \sort{
    \field{sorttitle}
    \field{title}
  }
}

\usepackage[utf8]{inputenc}
\usepackage[shortlabels]{enumitem}
\usepackage[final]{microtype}
\usepackage{bm}
\usepackage{csquotes}
\usepackage[english]{babel}
\usepackage{bbm}

\numberwithin{equation}{section}
\newtheorem{teo}{Theorem}[section]

\newtheorem{prop}[teo]{Proposition}
\newtheorem{lemma}[teo]{Lemma}
\newtheorem{cor}[teo]{Corollary}

\theoremstyle{definition}
\begingroup
\newtheorem{remark}[teo]{Remark}

\endgroup

\theoremstyle{remark}

\newcommand{\z}{\mathbb{Z}}

\newcommand{\restr}[2]{\left.\kern-\nulldelimiterspace#1\vphantom{\big|}\right|_{#2}}

\newcommand{\argmin}{ \mathrm{argmin} \,}
\newcommand{\F}{\mathcal F}
\newcommand{\G}{\mathcal G}
\newcommand{\K}{\mathcal K}
\newcommand{\strain}{\textbf{$\epsilon$}}

\DeclareMathOperator{\BV}{BV}
\DeclareMathOperator{\SBV}{SBV}

\Crefname{lemma}{lemma}{lemmas}
\Crefname{lemma}{Lemma}{Lemmas}
\Crefname{prop}{proposition}{propositions}
\Crefname{prop}{Proposition}{Propositions}
\Crefname{defi}{definition}{definitions}
\Crefname{defi}{Definition}{Definitions}
\Crefname{cor}{corollary}{corollaries}
\Crefname{cor}{Corollary}{Corollaries}
\Crefname{remark}{remark}{remarks}
\Crefname{remark}{Remark}{Remarks}
\Crefname{teo}{theorem}{theorems}
\Crefname{teo}{Theorem}{Theorems}
\Crefname{section}{Section}{Sections}
\Crefname{app}{Appendix}{Appendices}
\Crefname{example}{Example}{Examples}
 
\newcommand{\bl}{\color{black}}

\title[$\boldsymbol{\Gamma}$-convergence for a phase-field cohesive energy]{$\boldsymbol{\Gamma}$-convergence for a phase-field cohesive energy}
\author[E. Maggiorelli]{Eleonora Maggiorelli\hspace{.5mm}\orcidlink{0009-0009-7986-2087}}
\author[M. Negri]{Matteo Negri\hspace{.5mm}\orcidlink{0000-0001-6161-4020}}
\author[F. Vicentini]{Francesco Vicentini}\author[L. De Lorenzis]{Laura De Lorenzis}
\begin{document}

\begin{abstract}
Reproducing the key features of fracture behavior under multiaxial stress states 
is essential for accurate modeling. Experimental evidence indicates that three intrinsic material properties govern fracture nucleation in elastic 
materials: elasticity, strength, and fracture toughness (or critical energy release rate). Among these, strength remains the most often misunderstood, as it is not a single scalar quantity but rather a full surface in stress space. The flexibility in defining this strength envelope in phase-field models poses significant challenges, especially under complex loading conditions. 

Existing models in the literature often fail to capture both the qualitative shape and the quantitative fit of experimentally observed strength surfaces. To address this limitation, recent work introduces a new energy functional within a cohesive phase-field framework, specifically designed to control the shape of elastic domains. This model introduces an internal  variable to describe the inelastic response. Notably, the strength is decoupled from the internal length, that is not interpreted as a material length scale, as often done in literature, but rather as a purely variational tool.
The proposed functional allows for a rigorous variational framework, enabling the use of tools from the calculus of variations. We investigate the $\Gamma$-convergence of the model to a sharp cohesive fracture energy in the one- and two-dimensional (anti-plane) setting, using a finite element discrete formulation and exploiting the strong localization of the damage variable. Notably, unlike classical models where the elastic and fracture energies converge independently, this model exhibits a coupling of all energy terms.
The limiting cohesive energy arises from the combined asymptotic behavior of the elastic energy (concentrated in a single element), the fracture energy, and the  potential for the internal variable, while the remaining elastic energy converges separately.

We also present numerical simulations exploring the sensitivity of the model to mesh anisotropy, offering insight into both its theoretical robustness and its practical implementation.	
	\end{abstract}
	
	\maketitle
	
	{\small
		\par

		}

	\pagenumbering{arabic}
	
	\medskip
	
	\tableofcontents

\input{introduction}
\input{formulation}
\input{appendix}

	\bigskip
	
	\noindent\textbf{Acknowledgements.} E. Maggiorelli and M. Negri are members of GNAMPA - INdAM and are supported by PRIN 2022 (Project no. 2022J4FYNJ), funded by MUR, Italy, and the European Union -- Next Generation EU, Mission 4 Component 1 CUP F53D23002760006. F. Vicentini and L. De Lorenzis acknowledge funding from the Swiss National Science Foundation through Grant N. 200021-219407 ‘Phase-field modeling of fracture and
fatigue: from rigorous theory to fast predictive simulations’.

	\bigskip
    \printbibliography

	\medskip
	\small
	\begin{flushleft}
		\noindent \verb"eleonora.maggiorelli01@universitadipavia.it"\\
		Dipartimento di Matematica  ``F. Casorati'', Universit\`a di Pavia,\\
		via Ferrata 5, 27100 Pavia, Italy\\
		\medskip
		\noindent \verb"matteo.negri@unipv.it"\\
		Dipartimento di Matematica  ``F. Casorati'', 
		Universit\`a di Pavia,\\
		via Ferrata 5, 27100 Pavia, Italy\\
        \medskip
        \noindent \verb"fvicentini@ethz.ch"\\
		Computational Mechanics Group, Eidgen\"{o}ssische Technische Hochschule Z\"{u}rich,\\ 
        Tannenstrasse 3, 8092 Z\"{u}rich, Switzerland\\
        \medskip
        \noindent \verb"ldelorenzis@ethz.ch"\\
		Computational Mechanics Group, Eidgen\"{o}ssische Technische Hochschule Z\"{u}rich,\\ 
        Tannenstrasse 3, 8092 Z\"{u}rich, Switzerland\\
		\smallskip
	\end{flushleft}

\end{document}

%% file: introduction.tex
\section{Introduction}

Accurate modeling of fracture in brittle and quasi-brittle materials requires capturing both the nucleation of new cracks and the propagation of existing cracks, which are governed by two independent material properties: strength and toughness, respectively. Under multiaxial stress states, strength is not a single scalar quantity but is instead represented by a (convex) surface in stress space, separating admissible from inadmissible stress states. These strength surfaces are typically asymmetric, reflecting the markedly different behavior of materials in tension versus compression. Classical strength criteria that define such surfaces include those of Rankine, Mohr–Coulomb, and Drucker-Prager.

Griffith’s seminal fracture theory assumes the presence of an existing crack and—under certain simplifying assumptions, such as idealized geometries and planar crack paths—predicts whether the crack will propagate by comparing the energy release rate to the material’s fracture toughness. However, it does not account for the nucleation of new cracks. In Griffith’s framework, the energetic cost of fracture is assumed to be independent of the displacement jump across the crack faces and proportional only to the crack surface area. The variational reformulation of Griffith’s theory \cite{francfort1998} overcomes (at least at the theoretical level) the limitation of prescribing the crack path but still cannot adequately model crack nucleation.

For many years, strength-based criteria and fracture mechanics evolved as conceptually separate approaches. This gap was bridged with the advent of cohesive zone models (notably by Barenblatt \cite{barenblatt1962}, Dugdale \cite{dugdale1960}, and later Hillerborg \cite{Hillerborg1976}), which assign an energy cost to fracture that depends on the magnitude of the displacement jump across the crack faces. This framework enables a smooth transition from intact material to fully developed cracks, effectively unifying the modeling of crack initiation and propagation.   These cohesive models have since been the focus of extensive mathematical analysis, covering properties of minimizers \cite{DalMaso2008,Caffarelli2020}, evolutions in a one-dimensional setting \cite{Braides1999,Artina,Bonacini2021} and evolutions in the plane-strain setting along prescribed interfaces \cite{Almi2017,Negri2017}.

Both brittle and cohesive formulations lead to a challenging free-discontinuity problem, which is difficult to tackle numerically. Thus, implementing such problems requires suitable approximations of the energy, in the sense of $\Gamma$-convergence \cite{DalMaso93,Braides98}. In the case of brittle fracture, the Ambrosio-Tortorelli regularization \cite{BourdinFrancMar} - also interpretable as a damage model \cite{Pham2010} - provides an effective approximation by means of a separately convex energy, for which $\Gamma$-convergence has been proven in several settings  \cite{AmbrosioTortorelli1990,Bellettini1994,Chambolle2004,Chambolle2018}. As a by-product, global minimizers of the energy converge, as well as quasi-static evolutions \cite{FrancLars03,Giacomini2005,Giacomini2006}. However, global minima do not provide in general physically sound evolutions and numerical methods rather compute critical points or local minima. Convergence of critical points has been recently proven in \cite{BabadjianMillotRodiac_24}. A study of the evolutions in terms of critical points shows theoretically \cite{MaggiorelliNegri} and numerically \cite{Maggiorelli2025} that phase-field evolutions approximate sharp crack evolutions governed by Griffith's criterion and maximal energy release rate. In this perspective, phase-field approximations completely solve the problem of crack path selection.

Phase-field approximations also introduce nucleation (interpreted as the loss of second-order stability of nearly uniform damage solutions under local minimization \cite{TanLiBou17}). However, the resulting strength surface is elliptic, allowing only a single strength parameter—typically the tensile strength—to be calibrated through the choice of the regularization length, which effectively becomes a material parameter \cite{delor2022}. This is inadequate to capture the asymmetric tensile–compressive behavior observed experimentally. Furthermore, the model does not account for unilateral contact at the crack faces. 
To approximate more realistic strength surfaces and incorporate unilateral contact, several extensions based on energy decomposition have been proposed. The volumetric–deviatoric decomposition \cite{LancioniRoyerCarfagni2009,amor2009} admits a $\Gamma$-convergence result to brittle fracture with unilateral contact \cite{Chambolle2018}. Other decompositions  have been developed \cite{Miehe2010,freddi2010,delor2022,vicentini2024energy}, for which  $\Gamma$-convergence results are lacking, however the flexibility of these approaches remains limited. Among them, the recent model in \cite{vicentini2024energy} enables the separate calibration of tensile and compressive strength while satisfying all additional desirable requirements. However, it is limited to a strength surface of star-convex shape. Moreover, the strength surfaces obtained by all these models become unbounded as the regularization length tends to zero. The limited flexibility achievable by the energy decomposition approach is perhaps unsurprising, given that all these models approximate Griffith’s theory, which inherently lacks a notion of strength.

More recently, some regularized cohesive fracture models have been proposed \cite{CFI,Conti2024}. However, these models typically require modifications for numerical implementation \cite{FI}, and their mathematical structure remains very similar to that of the Ambrosio-Tortorelli functional. As a result, it is questionable whether the desired level of flexibility can be achieved within this framework.

To address the previous limitations, recent works \cite{vicent,Bourdin2025} introduce a new energy functional within a cohesive phase-field framework, specifically designed to control the shape of the strength surface and inspired by a similar functional introduced in \cite{AlessiMarigoVidoli_ARMA14}, for which 
$\Gamma$-convergence to a cohesive fracture model for an elasto-plastic material is shown in \cite{dalmaso2016}. These models introduce an internal  variable to describe the inelastic response, which is interpreted as an inelastic strain (for fracture in elasto-plastic materials such as in \cite{AlessiMarigoVidoli_ARMA14}, it coincides with the plastic strain). Notably, as a result of the new formulation, the strength  is decoupled from the regularization parameter, that is no longer a material property. In this paper we focus on this new energy functional. 

Besides the phase-field approach, it is worth to mention eigen-fracture \cite{SchmidtFraternaliOrtiz_MMS08,PO_IJNME12,ABS_AMPA22} which features a plastic-like variable, as our phase-field functional. Considering in particular finite element approximations for both eigen-fracture and our phase-field approach, cracks are represented (geometrically) by a narrow stripe of elements with large displacement gradient, while the fracture energy is computed by non-local terms, which prevent mesh bias. In phase-field the non-locality is obtained by setting the internal length to be (much) larger than the mesh size, so that the finite element solution can represent accurately the transition profile of damage. In eigen-fracture non-locality is introduced by means of non-local averages in a neighborhood (of the fracture elements) whose size plays the role of the internal length.

In this paper, we investigate the $\Gamma$-convergence of the functional in \cite{vicent,Bourdin2025} to a sharp cohesive fracture energy in the one- and two-dimensional (antiplane) settings, using a finite element discrete formulation and exploiting the strong localization of the damage variable.
The paper is organized as follows. Sections \ref{sec:formulation}  and \ref{sec:3} present the main $\Gamma$-convergence results in the one-dimensional and two-dimensional anti-plane settings, respectively. Section \ref{s.bib} provides an brief overview of related $\Gamma$-convergence results and compares them with the proposed model. In particular, we consider the eigen-fracture approximation of brittle fracture energies, introduced in \cite{SchmidtFraternaliOrtiz_MMS08}, as well as the eigen-fracture approximations of cohesive fracture energies, namely \cite{DMOT} and \cite{ABS_AMPA22}. We also mention the phase-field formulation developed in \cite{CFI,FI}, which, unlike our approach, does not rely on any additional internal variable.
Sections \ref{optprof}, \ref{sec:limsup}, and \ref{sec:liminf} are dedicated to the detailed proof of $\Gamma$-convergence in the one-dimensional setting. Finally, Section \ref{s.num} presents numerical simulations that validate the isotropy of the discrete energy. In particular, we investigate whether the finite element discretization introduces mesh-induced anisotropy, and show that the formulation remains robust with respect to mesh geometry.

%% file: formulation.tex
\section{$\Gamma$-convergence for the 1D model}\label{sec:formulation} 
We consider a mesh $\mathcal T_h$ of size $h$ in the domain $I=[-L,L]$ and denote $\mathbb{P}_h^0$ and $\mathbb{P}_h^1$ the spaces of piecewise constant and piecewise affine functions, respectively. Set $\epsilon_h >0$ with $h = o (\epsilon_h)$, 
we define the discrete functional
$F_h:\mathbb{P}_h^1\times \mathbb{P}_h^0 \times \mathbb{P}_h^1\to [0,+\infty)$ as follows:
$$F_h(u_h,\eta_h,d_h)=\int_I \tfrac12E_0|u_h'-\eta_h|^2 \, \text{d}x + \int_I a(\bar d_h)\sigma_c\eta_h\,\text{d}x+\tfrac{G_c}2\int_{I}\frac{1}{\epsilon_h}d_h^2+\epsilon_h|d_h'|^2\,\text{d}x .$$
The variables $u_h$ and $d_h$ are respectively the displacement variable and the damage variable, while $\eta_h \ge 0$ is used to introduce a threshold for the inelastic behavior. The notation $\bar d_h$ indicates the mean value of $d_h$ on each element. The degradation function is $a(d)=(1-d)^2$. 

\medskip
In the sequel it will be convenient to define 
the functionals $\F_h:\,\mathbb{P}_h^1\times \mathbb{P}_h^1\to[0,+\infty)$, depending only on the displacement $u_h$ and on the damage variable $d_h$:
\begin{equation}   \label{F_h}\F_h(u_h,d_h)=\min\{F_h(u_h,\eta_h,d_h)\,:\,\eta_h\in  \mathbb{P}_h^0,\,\,\eta_h\geq 0\}.
\end{equation}
Since $\eta_h\in \mathbb{P}_h^0$, this minimization can be done element by element. 
The functional $\F_h$ can hence be written in integral form as:
\begin{equation}
   \label{eq: Fhint}\F_h(u_h, d_h)=\int_I f(u_h',\bar d_h)\text{d}x+\tfrac{G_c}2\int_{I}\frac{1}{\epsilon_h}d_h^2+\epsilon_h|d_h'|^2\,\text{d}x 
\end{equation}
where 
\begin{align} \label{eq:f}
 f(s,r)&=\min\bigg\{\tfrac12 E_0 (s-\eta)^2+a(r)\sigma_c\eta \,\,:\,\eta \geq 0\bigg\}\\&=\begin{cases} \frac{1}{2} E_0 s^2 & s\leq a(r)\frac{\sigma_c}{E_0}, \\
a(r)\sigma_cs-\frac{\sigma_c^2}{2E_0}a^2(r) & s>a(r)\frac{\sigma_c}{E_0}.
\end{cases}  
\end{align}

On every element, the energy density $f ( \cdot , \bar d_h)$ is quadratic up to the threshold value $a(\bar d_h)\frac{\sigma_c}{E_0}$, that decreases as the damage increases.
In particular, where $\bar d_h=1$, it takes the form
$$f(u_h',1)=\begin{cases}
    \frac12E_0|u_h'|^2 \qquad\qquad \text{ for $u'_h\leq 0$}, \\ 0\qquad \qquad \quad \qquad\,\,\text{ for $u'_h\geq 0$}, 
\end{cases}$$ accounting for the loss of tensile strength in areas where the material is broken, while in compression the material is still elastic. On the contrary, where $\bar{d}_h = 0$, the elastic energy reads 
\begin{align} \label{eq:f}
  f(u'_h,0) & =\begin{cases} \frac{1}{2} E_0 | u'_h |^2 &  \text{ for $u'_h \le \frac{\sigma_c}{E_0}$}, \\
\sigma_c u'_h - \frac{\sigma_c^2}{2E_0} & \text{ for $u'_h \ge \frac{\sigma_c}{E_0}$} , 
\end{cases}  
\end{align}
showing a plastic-like behaviour under tension and a purely elastic behavior in compression. 

\medskip

As a boundary condition, we impose $u_h=g$ and  $d_h=0$ in $\partial I = \{ \pm L \}$. The latter ensures that, in the discrete setting, damage does not occur in the presence of a Dirichlet boundary condition for the displacement. Moreover, we require $d_h \in [0,1]$ and, without loss of generality, that $\|u_h\|_\infty\leq \| g \|_\infty  $.

\medskip
The asymptotic behavior of the functionals $\F_h$ is obtained by studying their $\Gamma$-limit in the space $L^1(I)\times L^1(I)$ as $h \to 0$.
The functionals \eqref{F_h} are thus extended to the space $L^1(I)\times L^1(I)$ by setting: $$\widetilde\F_h(u_h,d_h)=\begin{cases}\F_h(u_h,d_h) & \text{if $u_h,\,d_h\in \mathbb{P}_h^1$,\, $\| u_h \|_\infty \le \| g \|_\infty$,\,$d_h\in[0,1]$} \\ &\text{$u_h=g$, $d_h=0$ in $\partial I$} \\[2pt]+\infty & \text{otherwise.}
    \end{cases}$$ 
For $\Gamma$-convergence to hold, we require the mesh size to be sufficiently smaller than the internal length, i.e., $h=o(\epsilon_h) $; this ensures an accurate approximation of the transition layer of the phase-field variable and, in practice, it only needs to be satisfied in a neighbourhood of the discontinuity set. This is typically achieved through a local $h$-refinement. The theorem presented below constitutes the main result of this work and will be proved in \S\ref{sec:limsup} and \S\ref{sec:liminf}.

\begin{teo}\label{TH} As $h \to 0$, the functionals $\widetilde\F_h$ $\Gamma$-converge to $\widetilde{\F}: L^1(I)\times L^1(I)\to [0,+\infty]$  defined as follows:
\begin{equation}
  \label{eq: F}\widetilde{\F}(u,d)=\begin{cases}
    \F(u) & \text{if } u\in \BV(I),\, \| u \|_\infty \le \| g \|_\infty,\,\llbracket u\rrbracket>0,\,D^cu\geq 0, \text{ and }d=0 \text{ a.e.~in $I$,}\\
    +\infty & \text{otherwise,}
\end{cases}    
\end{equation} 
where 
$$\F(u)=\int_I W(u')\, \text{d}x+\sigma_c|D^cu|+\sum_{J_u}\phi(\llbracket u\rrbracket )+\phi(g(L)-u(L^-))+\phi(u(-L^+)-g(-L))$$
and  the functions $W$ and $\phi$ are given by:
   \begin{equation}
       \label{eq:phi} W(s) = f ( s, 0), \quad \phi(s)= \begin{cases}
\displaystyle          \frac{G_c\sigma_cs}{G_c+\sigma_cs} & \text{if $s\ge0$} \\+\infty & \text{if $s<0$.}
       \end{cases}
   \end{equation}
   \end{teo}

\begin{remark} Note that for $s \ge 0$ the cohesive energy $\phi$ is concave and increasing, with $\phi'_+ (0) = \sigma_c$ and $\lim_{s \to +\infty} \phi(s) = G_c$. 
Even if the discrete damage variable $d_h$ is null in $\partial I$, in the limit, damage at the boundary can still occur, i.e., we may have $d \neq 0$ in $\partial I$.  Indeed, in $L^1(I)$,  we can approximate with finite energy a function that is non-zero in $L$ (or, equivalently, $-L$) using functions that vanish at the boundary. As a consequence we may as well have $u \neq g$ in $\partial I$. The definition of $\F$ above highlights that we are considering all jump points, including those at the boundary. However, for the sake of readability, it is convenient to express the functional in a more compact form. To do so, we introduce the following function:
$$\label{utilde}\tilde u(x)=\begin{cases}g(L) \qquad x\geq L\\
u(x)\qquad x\in(-L,L),\\
g(-L)\quad\, x\leq-L. 
\end{cases}$$
This allow us to define:
\begin{equation}
    \F(u)=\int_I f(u',0)\, \text{d}x+\sigma_c|D^cu|+\sum_{J_{\tilde u}}\phi(\llbracket\tilde u\rrbracket).
\end{equation}
By employing this extension of $u$ while setting homogeneous boundary conditions on the damage variable, we ensure that the fracture energy contribution at the boundary is not artificially reduced.
\end{remark}

\section{$\Gamma$-convergence for the 2D antiplane model}\label{sec:3}

In this section we state the $\Gamma$-convergence result in the anti-plane case, as in \cite{DMOT}. For the sake of simplicity, let $\Omega = (-L,L) \times (-H,H)$ and let $\mathcal T_h$ be a regular triangulation of the domain. By abuse of notation, we still denote by $\mathbb{P}_h^i$ ($i=0,1$) the spaces of piecewise constant and piecewise affine finite elements on $\mathcal T_h$. 

The discrete energy $F_h : \mathbb{P}_h^1 \times (\mathbb{P}_h^0 \times \mathbb{P}_h^0  ) \times \mathbb{P}_h^1 \to \mathbb{R}$ is then given by 
\begin{equation} \label{e.defFh}
   F_h ( u_h , \eta_h , d_h) = \int_\Omega \mu | \nabla u_h - \eta_h |^2 \, \text{d}x + 
   \int_\Omega a ( \bar{d}_h ) \sigma_c | \eta_h | \, \text{d}x + 
   \tfrac{G_c}{2} \int_\Omega \frac{d_h^2}{\epsilon_h} + \epsilon_h | \nabla d_h |^2 \, \text{d}x ,
\end{equation}
where $\mu>0$ is the shear modulus and $\epsilon_h>0$ with $h = o (\epsilon_h)$. Being $\eta_h$ a piecewise constant vector field, it is convenient, as in \S \ref{sec:formulation}, to minimize the energy density on every element, and then consider $\mathcal{F}_h : \mathbb{P}_h^1 \times \mathbb{P}_h^1 \to [0,+\infty)$ given by 
$$
\mathcal{F}_h (u_h, d_h) = \int_\Omega f ( | \nabla u_h | , d_h) \, \text{d}x + \tfrac{G_c}{2} \int_\Omega \frac{d_h^2}{\epsilon_h} + \epsilon_h | \nabla d_h |^2 \, \text{d}x ,
$$
where, in analogy with \eqref{eq:f}, by radial symmetry $f$ is given by 
\begin{align} \label{eq:fbis}
 f(s,r)&=\min \big\{\mu (s-\eta)^2+a(r)\sigma_c \eta  \,\,:\,\eta \ge 0 \big\}\\&=\begin{cases} \mu s^2 & s\leq a(r)\frac{\sigma_c}{2 \mu }, \\
a(r)\sigma_cs-\frac{\sigma_c^2}{4\mu}a^2(r) & s>a(r)\frac{\sigma_c}{2 \mu }.
\end{cases}  
\end{align}


\medskip
For $g \in H^1(\Omega) \cap L^\infty(\Omega)$ we consider the Dirichlet boundary conditions $u_h = g$ and $d_h = 0$ in $\partial_D \Omega = \{ \pm L \} \times (-H, H)$. Moreover, we consider the constraint $\| u_h \|_\infty \le \| g \|_\infty$ and $d_h \in [0,1]$. Then, the extended functional $\widetilde{\mathcal{F}}_h : L^1 (\Omega) \times L^1(\Omega) \to [0,+\infty]$ is given by
$$\widetilde\F_h(u_h,d_h)=\begin{cases}\F_h(u_h,d_h) & \text{if $u_h,\,d_h\in \mathbb{P}_h^1$,\, $\| u_h \|_\infty \le \| g \|_\infty$,\, $d_h\in[0,1]$} \\ &\text{$u_h=g$, $d_h=0$ in $\partial_D \Omega$} \\[2pt]+\infty & \text{otherwise.}
    \end{cases}$$ 
At this point, we can state the $\Gamma$-convergence result, considering again $h = o (\epsilon_h)$. 

\begin{teo}\label{TH2D} As $h \to 0$, the functionals $\widetilde\F_h$ $\Gamma$-converge to $\widetilde{\F}: L^1(I)\times L^1(I)\to [0,+\infty]$  defined as follows:
\begin{equation}
  \label{eq: F}\widetilde{\F}(u,d)=\begin{cases}
    \F(u) & \text{if } u\in \BV(\Omega),\, \| u \|_\infty \le \| g \|_\infty,\,\text{ and $d=0$ a.e.~in $\Omega$,}\\
    +\infty & \text{otherwise,}
\end{cases}    
\end{equation} 
where 
$$\F(u)=\int_\Omega f(| \nabla u| ,0)\, \text{d}x+\sigma_c|D^c u|+ \int_{J_u} \phi (| \llbracket u\rrbracket | )+ \int_{\partial_D \Omega} \phi(| u - g | )$$
and $\phi$ is defined in \eqref{eq:phi}. 
   \end{teo}

\begin{remark} \label{iso}In this setting the non-interpenetration condition does not apply and indeed both positive and negative jumps are allowed. In the plane-strain setting, the non-interpenetration condition is instead a difficult technical point, preventing a complete $\Gamma$-convergence result.
\end{remark}

\begin{remark} The limit energy \eqref{eq: F} is isotropic, i.e., it is independent of the geometry of the underlying triangulation. This property is confirmed in the numerical simulations of \S \ref{s.num}, actually performed in plane strain. Noteworthy, in accordance with our $\Gamma$-convergence proof, displacement jumps are approximated at the element size, while the fracture energy depends on the phase-field profile in a neighborhood of size $\epsilon_h \gg h$. This prevents the mesh bias in analogy with non-local averaging in eigen-fracture \cite{SchmidtFraternaliOrtiz_MMS08} and smeared crack \cite{LussN_NFAO07} approaches.
\end{remark}

\section{Some related $\Gamma$-convergence results\label{s.bib}}

In this section we briefly discuss the relationship between our result and: (a) the approximation \cite{SchmidtFraternaliOrtiz_MMS08} of brittle energies and (b) a couple of approximations of cohesive energies, specifically \cite{DMOT} and \cite{ABS_AMPA22}. All these results share the use of a ``plastic-like variable'' but they also have interesting differences. We finally mention the phase-field approximation \cite{CFI,FI} which however does not employ a plastic variable.

\newcommand{\eps}{\epsilon}

Let us start from the {\it eigen-fracture} approach of \cite{SchmidtFraternaliOrtiz_MMS08}. To better compare with our result, we consider the anti-plane finite element discretization, which is enough to show that mesh bias does not occur (see also \S \ref{s.num}). Let $h = o (\epsilon_h)$, as in our setting. Given $A \subset \Omega$ let $A_h$ denote the union of the elements $e_h \in \mathcal{T}_h$  such that $\mathrm{dist} ( e_h , A ) \le \epsilon_h$. 
 In our notation, the functional $F_h :  \mathbb{P}^1_h \times ( \mathbb{P}^0_h \times \mathbb{P}^0_h) \bl \to [0, +\infty]$ introduced in \cite{SchmidtFraternaliOrtiz_MMS08} takes the form 
\begin{equation} \label{e.FOS}
	F_h ( u_h , \eta_h) =  \int_\Omega \mu | \nabla u_h - \eta_h |^2 \, \text{d}x 
		+ \tfrac{G_c}{2 \epsilon_h} \, \big|  \{  \eta_h \neq 0   \}_h  \big| ,
\end{equation}
where for simplicity we neglect here the boundary condition and the bound on $\| u_ h \|_\infty$. Note that the measure term $| \{ \eta_h \neq 0 \}|$ is not differentiable with respect to $\eta_h$.  In this approximation the ``plastic variable'' $\eta_h$ is again concentrated (see \cite{SchmidtFraternaliOrtiz_MMS08,PO_IJNME12}) on a single stripe of elements (of order $h$) while the surface energy depends on the internal lenght $\eps_h$ as in the phase-field approach. This non-locality allows to avoid mesh bias in the approximation of the surface energy.  Indeed, for $h = o (\eps_h)$ the energy  $F_h$ $\Gamma$-converges (as $h \to 0$) to the Griffith energy
$$ 
	\mathcal{F} (u) = \int_\Omega \mu |  \nabla u |^2  + G_c \mathcal{H}^1 ( J_u) .
$$

Next, let us consider \cite{DMOT} and \cite{ABS_AMPA22}. Note that both results are set in the spatially continuum setting and there is no finite element discretization. Moreover they do not consider the unilateral constraint on the crack.  As a common root, we restrict to the one dimensional formulation with quadratic fracture energy, which is however enough to characterize the cohesive energy density.

The convergence result of \cite{DMOT} provides a rigorous mathematical proof of the {\it phase-field} energy originally proposed in \cite{AlessiMarigoVidoli_ARMA14}. The phase-field energy $F_\eps : BV(I) \times \mathcal{M} (I) \times H^1(I; [0,1] ) \to [0,+\infty]$ takes the form 
\begin{equation} \label{e.DOT}
	F_\eps ( u , \eta , d ) = \int_I (1-d) | e |^2 \, \text{d}x + \int_I  k (d) \text{d} | \eta |  + \tfrac{G_c}{2} \int_I \tfrac{1}{\eps} d^2 + \eps | d' |^2 \, \text{d}x, 
\end{equation}
where $u' = e + \eta$ with $e \in L^2$. 
Comparing with \eqref{e.defFh} note that here the displacement $u$ is not defined in the Sobolev space $H^1$ but in the larger space $BV$, moreover, the elastic energy features the degradation function $(1-d)$, while $k(d)$ plays the role of $a(d) \sigma_c$. Minimizing with respect to $\eta$ provides the reduced functional 
$$
	\mathcal{F}_\epsilon ( u , d) = \int_I f_\eps (u', d) \, \text{d}x + \int_I  k (d) \, \text{d} | D^s u |  + \tfrac{G_c}{2}  \int_I \tfrac{1}{\eps} d^2 + \eps | d' |^2 \, \text{\text{d}x}, 
$$
where $Du = u' + D^s u$ while $f_\eps$ has a quadratic-linear structure similar to \eqref{eq:f}. 
Choosing $k(d) = \sigma_c (1-d)^2$ the $\Gamma$-limit (as $\eps \to 0$) is given by the functional 
$$
	\mathcal{F} (u) = \int_I  f ( u'  , 0)  \, \text{d}x + \sigma_c | D^c u | + \sum_{x \in J_u} \phi \big( \llbracket u \rrbracket  \big) , 
$$
where $\phi (s) = G_c \sigma_c |s| / (G_c + \sigma_c |s|) $ coincides with \eqref{eq:phi} for $s \ge 0$, while $f$ coincides with $W$ for $E_0=1$. Comparing with our result it turns out that in the discrete setting it is not restrictive to consider displacements $u_h$ in $H^1$, instead of $BV$, and that it is not necessary to employ degradation functions, as in eigen-fracture models.  Our convergence proof is crafted for the discrete setting and is indeed independent of that of \cite{DMOT}.

The energy studied in \cite{ABS_AMPA22} has its root in the {\it eigen-fracture} approach \cite{SchmidtFraternaliOrtiz_MMS08} described above and in the non-local approximation of \cite{LussardiVitali_AMPA07}.  In this case the energy $F_\epsilon: BV(I) \times {\mathcal M} (I) \to [0,+\infty]$ is given by 
$$
	F_\eps  (u, \eta) = \int_I \tfrac{E_0}2  | u' - \eta |^2 \, \text{d}x + 
	\tfrac{1}{2 \epsilon} \int_I \varphi \Big( \int_{(x-\eps,x+\eps)} | \eta | \, \text{dy} \Big) \, \text{d}x ,
$$
where $\eta \in L^1$ with $(u' - \eta) \in L^2$ while $ \varphi(s) = \sigma_c  
\min \{  |s| , 1 \}$. 
Once again, minimizing with respect to $\eta$, the $\Gamma$-limit (as $\eps \to 0$) takes the form
$$
	\mathcal{F} (u) = \int_I f ( u' , 0 ) \, \text{d}x +  \sigma_c | D^c u | +\sum_{x \in J_u}  \varphi \left( \llbracket u \rrbracket \right) . 
$$

For the sake of completeness we mention also the phase-field approximation \cite{CFI,FI} which does not employ any plastic variable but a suitable degradation function, i.e.
$$
	F_\epsilon ( u, d) = \int_I f_\eps ( u' , d)\, \text{d}x + \tfrac{G_c}{2} \int_I \tfrac{1}{\eps} d^2 + \eps | d' |^2 \, \text{d}x . 
$$ 
Here, the function $f_\eps$ takes the form $f_\eps (s,r) = |s|^2 \min \{ \eps^{1/2} \psi (r) , 1 \}$ with $\lim_{r \to 0^+} r \psi (r) = \sigma_c$.  Note that $\psi_\eps$ in general is non-convex and depends on the internal length $\epsilon$. The $\Gamma$-limit (as $\epsilon \to 0$) is a cohesive energy the form 
$$
	\mathcal{F} (u) = \int_I f ( u' ) \, \text{d}x + \sigma_c | D_c u | + \sum_{x \in J_u} \phi ( \llbracket u \rrbracket) 
$$
where $f$ has a quadratic-linear behaviour while the cohesive potential $\phi$ can be characterized in terms of $\psi_\eps$, appearing in $f_\eps$. 


\section{Optimal profile 
} \label{optprof}

\begin{figure}[h]
\centering
\includegraphics[scale=1]{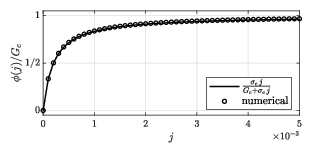}\caption{Analytical vs. numerical surface energy density $\phi$.}\label{fig:phi_1D}\end{figure}

Before proving the main convergence result, it is necessary to define the optimal profile problem. For the sake of readability, we start by defining the problem in a continuous setting and then consider its discrete approximation. 
Let us consider a solution of pure jump of positive amplitude, with $J_{u}=\{0\}$ and $\llbracket u \rrbracket(0)=j>0$. 
The optimal profile problem is the following: 
$$\label{z_j}z_{j}=\argmin\bigg\{J_{j}(z)=a(z(0))\sigma_cj+G_c\int_{\mathbb{R}^+}z^2+|z'|^2\,\text{d}x,\quad z \in H^1(\mathbb{R}_+, [0,1])\bigg \}.$$
To solve this problem, we first consider $z(0)=z_0$ as a fixed parameter and introduce the transition energy with unit internal length $\K: \mathcal{D}\to \mathbb{R}$
$$\K(z)=
\int_{\mathbb{R}^+}z^2+|z'|^2 \text{d}x,$$
where $\mathcal{D}=\{z\in H^1(\mathbb{R}_+,[0,1])\,:\,z(0)=z_0\}. $
The function $z_*(x)=z_0 e^{-x}$ is the unique minimizer of $\mathcal{K}$ over $\mathcal{D}$ 
and $\K(z_*)=
z_0^2$.
Therefore, to find $z_j(0)$, i.e the amplitude of the optimal profile, 
we need to solve:
$$z_{j}(0) \in \argmin \bigg\{ a( z_0)\sigma_cj+G_c z_0^2, \, \,\,z_0\in[0,1]\bigg\}.$$
By obvious calculations, 
$$z_{j}(0)=\frac{\sigma_c j}{G_c+\sigma_c j}\in[0,1]$$
and finally, we set $$\phi(j)=J_{j}(z_{j})=\frac{G_c\sigma_c j}{G_c+\sigma_c j}.$$ 
The surface energy density $\phi(j)$ is plotted in Fig.~\ref{fig:phi_1D}, showing the agreement between the analytical expression derived in this section (solid line) and the results from the numerical test described in Appendix~\ref{subsec:phi_num} (dots).

\section{Limsup-inequality } \label{sec:limsup}
\begin{prop} \label{limsup} Let $u\in \BV(I)$ such that $\| u \|_\infty \le \| g \|_\infty$, $\llbracket u \rrbracket>0$, and $D^cu \ge 0$. There exist $u_h, d_h \in \mathbb{P}_h^1\subset H^1(I)$ such that $(u_h,d_h)\to (u,0)$ as $h\to 0 $ in $L^1(I)\bl\times L^1(I)$ and:
    $$\limsup_{h \to 0}\F_{h}(u_h,d_h)\leq \F(u).$$ 

\end{prop}
\proof 
We recall the density result stated in \Cref{lem: denisty} and \Cref{rem-dens} which allows us to reduce the analysis to two representative cases: a pure jump function \( u \), and a smooth function \( u \in W^{2,\infty}(I) \).

{\bf I.} Let us consider the case in which $J_u=\{0\}$ and $u$ is constant elsewhere.
Assume without loss of generality that $u$ is left continuous. 
We set $u_h$ to be the piecewise affine interpolate of $u$, so that $u'_h=0$ on every element of the mesh except on the element $\bar e_h$ that contains the jump, where $u_h'=\frac{\llbracket u \rrbracket}{h}$. As shown in \Cref{optprof}, the solution of the optimal profile problem given $\llbracket u \rrbracket$ is $$z_{\llbracket u \rrbracket}(x)=\frac{\sigma_c \llbracket u \rrbracket}{G_c+\sigma_c \llbracket u \rrbracket}e^{-x}$$
and
$$\phi(\llbracket u \rrbracket)=J_{\llbracket u \rrbracket}(z_{\llbracket u \rrbracket})\leq J_{\llbracket u \rrbracket}(z)=a(z(0))\sigma_c\llbracket u \rrbracket+G_c\int_{\mathbb{R}^+}z^2+|z'|^2\,\text{d}x$$ 
for every $z\in H^1(\mathbb{R}_+,[0,1])$. For a fixed $\eta>0$, there exist $T_\eta>0$ and $z_\eta\in W^{1,\infty}(0,T_\eta)$
such that $z_\eta(0)=z_{\llbracket u \rrbracket}(0)$, \,$z_\eta(T_\eta)=0$, and
$$a(z_\eta(0))\sigma_c\llbracket u \rrbracket+G_c\int_{0}^{T_\eta}z_\eta^2+|z_\eta'|^2\,\text{d}x\leq \phi(\llbracket u \rrbracket)+\eta.$$
By abuse of notation, we call $z_\eta$ the null extension of such function to $T_\eta>0$.
Let us consider a rescaling of the mesh \( \mathcal T_h|_{[0,L]} \) by a factor \( 1/\epsilon_h \), and denote by \( \mathcal T_{h/\epsilon_h} \) the resulting mesh, now defined on the interval \( \left[0, \frac{L}{\epsilon_h}\right] \). Since $h=o(\epsilon_h)$, the mesh size $\frac{h}{\epsilon_h}\to 0$ as $h\to 0$ and for $h$ small enough $\frac{L}{\epsilon_h}>T_\eta$. 
We define \( z_h \) as the piecewise affine interpolate of \( z_\eta \) over the mesh \( \mathcal T_{h/\epsilon_h} \). Then, by standard finite element estimates, $z_h\to z_\eta$ strongly in $H^1(0,+\infty)$. Therefore for $h$ sufficiently small
\begin{equation} \label{e.kan} a\bigg(\frac{\sigma_c \llbracket u \rrbracket}{G_c+\sigma_c \llbracket u \rrbracket}\bigg)\sigma_c\llbracket u \rrbracket+G_c\int_{0}^{L/\epsilon_h}z_h^2+|z_h'|^2\,\text{d}x\leq \phi(\llbracket u \rrbracket)+2\eta.
\end{equation}
Finally, we introduce $d_h(x)=
    z_h( |x|/\epsilon_h)$. 
By definition $f(u_h',\bar d_h)=0$ outside $\bar e_h$, since $u'_h=0$. On $\bar e_h$ for $h$ sufficiently small we have instead \ $u_h'=\frac{\llbracket u \rrbracket}{h}>a (\bar d_h|_{\bar e_h})\frac{\sigma_c}{E_0}$, and hence $f|_{\bar e_h}$ is affine. In summary:
$$
f(u'_h,\bar d_h)=\begin{cases}
    0 \qquad \qquad \qquad\qquad \qquad \qquad\qquad I\setminus \bar e_h\\
a(\bar d_h)\sigma_c \frac{\llbracket u\rrbracket}{h}-\frac{\sigma_c^2}{2E_0}a^2(\bar d_h)\qquad\qquad\, \,\,\bar e_h.\end{cases}$$
By \eqref{eq: Fhint},
\begin{align}
    \F_h(u_h,d_h)&=hf\big(u'_h, \bar d_h|_{\bar e_h}\big)+\tfrac{G_c}2\int_{I}\frac{1}{\epsilon_h}d_h^2+\epsilon_h|d_h'|^2\,\text{d}x
    \\&=a(\bar d_h|_{\bar e_h})\sigma_c \llbracket u\rrbracket-h\frac{\sigma_c^2}{2E_0}a^2(\bar d_h|_{\bar e_h})+G_c\int_{0}^{L/\epsilon_h}z_h^2+|z_h'|^2\,\text{d}x .
    \end{align}
Note that 
$$
\bar d_h|_{\bar e_h} =  \tfrac12 ( z_\eta( 0) + z_\eta ( h / \epsilon_h) ) 
\to z_\eta(0) =
\frac{\sigma_c \llbracket u \rrbracket}{G_c+\sigma_c \llbracket u \rrbracket} . 
$$ 
Since $a$ is Lipschitz continuous on $[0,1]$, 
then for $h$ sufficiently small by \eqref{e.kan} it holds:
    \begin{align}
    \F_h(u_h,d_h)&\leq a\bigg(\frac{\sigma_c \llbracket u \rrbracket}{G_c+\sigma_c \llbracket u \rrbracket}\bigg)  \sigma_c \llbracket u \rrbracket 
+\eta+G_c\int_{0}^{L/\epsilon_h}z_h^2+|z_h'|^2\,\text{d}x\leq \phi(\llbracket u \rrbracket)+3\eta . 
\end{align}
 By the arbitrariness of $\eta$ we conclude that for a function of pure jump:
$$\limsup_{h\to 0}\F_h(u_h,d_h)\leq \phi(\llbracket u \rrbracket) = \F (u) .$$

{\bf II.} If $u\in W^{2,\infty}(I)$, we set $u_h$ to be the piecewise affine interpolate of $u$ and $d_h=0$. 
Then, by standard finite element estimates, see e.g. \cite{BS}, it holds 
$$\|u_h-u\|_{W^{1,\infty}}\leq h|u|_{W^{2,\infty}},$$
thus $u_h'\to u'$ uniformly and
$$\F_h(u_h,d_h)=\int_I f(u'_h,0)\text{d}x\to\int_I f(u',0)\text{d}x = \F(u),$$
which concludes the proof.  \qed
\begin{remark}   
    Notably, unlike classical models where the elastic and fracture energies converge independently, this model exhibits a coupling of all energy terms.
The limiting cohesive energy arises from the combined asymptotic behavior of the elastic energy (concentrated in a single element), the fracture energy and the plastic potential, while the remaining elastic energy converges separately. 
\end{remark}

\section{Liminf-inequality}
\label{sec:liminf}
\begin{prop} \label{liminf} 
    Let $(u_h,d_h)\to (u,d)$ as $h\to 0 $ in $L^1(I)\times L^1(I)$ and $\liminf_{h \to 0}\F_{h}(u_h,d_h)<+\infty$. 
    Then $d=0$ a.e. in $I$, $u\in \BV(I)$ with $\llbracket   u\rrbracket >0$ and 
    $$\F(u)\leq \liminf_{h \to 0}\F_{h}(u_h,d_h).$$
    As a direct consequence, $D^cu\geq 0$.
\end{prop}
\proof

The proof is carried out in several steps.

\textbf{I.} 
We start by proving the properties of the limit functions.
Since $(u_h,d_h)\to (u,d)$ as $h\to 0 $, there exist 
(non relabeled) subsequences of $(u_h, d_h)$ that converge to $(u,d)$ almost uniformly and such that $\lim_{h \to 0}\F_{h}(u_h,d_h)=\liminf_{h \to 0}\F_{h}(u_h,d_h)<+\infty$.
Therefore, from the fact that 
$$\frac{1}{\epsilon_h}\int_Id_h^2\text{d}x\leq \F_{h}(u_h,d_h)\leq C$$ 
and $\epsilon_h\to 0$, it follows that
$d=0$ a.e. in $I$. Since $\| u_h \|_\infty \le \| g \|_\infty$ it follows that $\| u \|_\infty \le \| g \|_\infty$. 
In the following, we prove that \( u \in \BV(I) \).  
For each \( m \in \mathbb{N} \), we introduce a set of ordered points $X^m=\{x_i^m\}_{i=0}^{m+1}$ with 
\( x^m_0 = -L \), \( x^m_{m+1} = L \) and such that: 
\begin{itemize}
    \item  \( d_h(x^m_i) \to 0 \) as \( h \to 0 \) for every \( i = 0, \dots, m+1 \); 
    \item \( \sup_i |I^m_i| \to 0 \) as \( m \to \infty \), where $| I^m_i|$ denotes the length of the sub-interval $I^m_i$.
\end{itemize}
 For every $i=1,...,m$, there exists an element in $\mathcal T_h$ which contains $x_i^m$. In this element, since $d_h\in \mathbb{P}_h^1$ in (at least) one of the endpoints, which we call $x^m_{i,h}$, we have $0 \le d_h(x^m_{i,h})\leq d_h(x^m_i)$. Moreover, we set $x_{0,h}^m=-L$ and $x_{m+1,h}^m=L$, and let  $I^{m}_{i,h}= 
[x^m_{i,h},x^m_{i+1,h}]$ for every $i\in \{0,...,m\}$. Clearly, $|x^m_{i,h}-x^m_i|<h$. 


 
First of all, note that
for every $m$ there exists $h_m$ such that $d_h(x^m_{i,h})\leq d_h(x^m_i)<1/4$ for every $i=0,..,m+1$ and every $h<h_m$. We now show the following: there exists $N>0$ such that for every $m\in\mathbb{N}$ and $h < h_m$ it holds: 
$$ N^m_{h} = \#\big\{I^{m}_{i,h}%
\,:\, \sup_{\,I^{m}_{i,h}} \{d_h\}> 1/2 \big\}< N . $$
%
We consider $h<h_m$ 
and we estimate the fracture energy on each interval $I_{i,h}^{m}$ on which $\sup \{ d_h \} > 1/2$. Since $\inf_{\,I^{m}_{i,h}}\{d_h\} < 1/4$, the fracture energy must exceed the minimal energy required to make a transition between $1/2$ and $1/4$. Therefore, recalling the optimal profile study of \Cref{optprof} we have:
$$ \tfrac{G_c}2\int_{I_{i,h}^{m}}\frac{1}{\epsilon_h}d_h^2+\epsilon_h|d_h'|^2\,\text{d}x\geq \min\{G_c\K(z)\,:\,z\in H^1(\mathbb{R}_+),\,z(0)=1/4 \}= c \, G_c .$$
As a consequence, for $h < h_m$:
$$c \,  G_c N_{h}^m \leq \F_h(u_h,d_h) \leq C$$
and thus we can set $N = C  / c \, G_c$. 
We call $I_{h}^m=\bigcup\{I^{m}_{i,h}
\,:\,\sup_{\,I^{m}_{i,h}}\{d_h\}> 1/2 \}$ and define
$$u^m_{h}=\begin{cases}
    u_h & I\setminus I_{h}^m%
    \\u_h(x^m_{i,h})+\big(u_h(x^m_{i+1,h})-u_h(x^m_{i,h})\big)\chi_{[\hat x_i,x^m_{i+1,h}]} &  I_{i,h}^{m}\subset I_{h}^m%
\end{cases}$$
for a some $\hat x_i\in (x^m_{i,h},x^m_{i+1,h})$.
We now show that $u^m_{h}$ is bounded in $\BV(I)$, studying separately the behavior in the subsets $I \setminus I_{h}^m$ and $I_{h}^m$.  

On $I\setminus I_{h}^m$, we have 
$\|d_h\|_\infty<1/2$ and hence, by convexity of $f( \cdot, 1/2 )$, we obtain 
$$ C \ge \int_{I\setminus I_{h}^m}f(u'_h,\bar d_h)\geq \int_{I\setminus I_{h}^m}f(u'_h,1/2)\geq \int_{I\setminus I_{h}^m}a(1/2)\sigma_c| u'_h| -\frac{\sigma_c^2}{2E_0}a^2(1/2)\,\text{d}x. $$
Then $u'_{h,m}=u'_h$ is bounded in $L^1({I\setminus I_{h}^m})$ and   
so  $|Du^m_{h}|(I\setminus I_{h}^m)$ is bounded.  
On the other hand, 
$$|Du^m_{h}|(I_{h}^m)
=\sum_{i=1}^{N_{h}^m}
|\llbracket u^m_{h}(\hat x_i)  \rrbracket|  \leq 2 N \|g\|_\infty.$$
As a consequence, $u^m_{h}$ is bounded in $\BV(I)$ for $h<h_m$ and up to non relabeled sequences, there exists a limit $u^m$ in $\BV(I)$. Since $u^m_{h}=u_h$  on $I\setminus I_{h}^m$ and by hypothesis $u_h\to u$ in $L^1(I)$, the limit $u^m$ must be equal to $u$ on $I\setminus \lim_{h\to 0} I_{h}^m=: I\setminus I^m$, where the set $I^m$ is the union of at most $N$ intervals $I^m_i$. Hence 
$ \|u\|_{\BV(I\setminus I^m)}<C $, where $C$ is independent of $m$. Now, since $\sup_i|I^{m}_i|\to 0$ as $m\to+\infty$, then $I^m\to\bigcup_{j=1}^{n}\{x_{j}\}$ where $n <N$ and hence 
$$\|u\|_{\BV(I\setminus\bigcup_{j=1}^{n}\{x_{j}\})}<C.$$
Finally, the finiteness of $\bigcup_{j=1}^{n}\{x_{j}\}$ and the fact that $\| u \|_\infty \le \| g \|_\infty$ ensures that $u\in \BV(I)$. 

It remains to show that $\llbracket u \rrbracket > 0$. We argue by contradiction, assuming  $\llbracket u (x) \rrbracket < 0$ for some $x \in I$. Let $\gamma>0$ (arbitrarily small) such that $u(x+\gamma) - u(x-\gamma) < 0$ and $u_h (x \pm \gamma) \to u (x \pm \gamma)$.  Denoting $I_\gamma = (x-\gamma, x+\gamma)$, we have
$$
    \int_{I_\gamma} f(u'_h , \bar{d}_h) \, \text{d}x \ge \inf \bigg\{ \int_{I_\gamma} f ( v',  \bar{d}_h ) \,\text{d}x \,:\,   v \in H^1(I_\gamma) \,,\, v (x\pm\gamma) = u_h(x\pm \gamma)  \bigg\} .
$$
For $h$ small enough $u_h(x+\gamma) - u_h(x-\gamma) < 0$, hence for minimality it is not restrictive to consider $v' \le 0$. It follows that
\begin{align} 
     \int_{I_\gamma} f(u'_h , \bar{d}_h) \, \text{d}x 
     & \ge \inf \bigg\{ \int_{I_\gamma} \tfrac12 E_0 | v'|^2 \, \text{d}x \,:\, v \in H^1(I_\gamma) \,,\, v' \le 0 \,,\, v(x\pm\gamma) = u_h(x\pm \gamma)  \bigg\} 
     \\ & \ge \tfrac{E_0}{4\gamma} | u_h (x+\gamma) - u_h (x-\gamma) |^2 .
\end{align}
The right hand side diverges as $\gamma \to 0$ which contradicts the boundedness of $\F_h (u_h, d_h)$. The fact that $D^c u \ge 0$ is not needed in the rest of the proof and it will follow from the liminf inequality itself. 

\medskip
\textbf{II.} We now prove that $$\F(u)\leq \liminf_{h \to 0}\F_{h}(u_h,d_h). $$

\textbf{Around jump points.} For $t>0$, we set $J^t_u=\{x\in J_u\,:\,\llbracket \tilde u\rrbracket\geq t\}$ and observe that, since $u\in\BV(I)$, $N^t=\#J_u^t<+\infty$. %
For $\delta >0$ sufficiently small the sets $J_u^{t,\delta}=\{x\in I\,:\, \mathrm{dist}(x,J_u^t)\leq \delta\}$ are disjoint intervals. For $I'\subset I$, we denote
$$\F_h(u_h,d_h,I')=\int_{I'} f(u_h',\bar d_h)\text{d}x+\tfrac{G_c}2\int_{I'}\frac{1}{\epsilon_h}d_h^2+\epsilon_h |d_h'|^2\,\text{d}x,$$
and since
$$\F_h(u_h,d_h)\geq \sum_{x \in J_u^t} \F_h(u_h,d_h,J_u^{t,\delta}),$$
we focus on a single $x_0\in J_u^t$.

If $x_0\neq \pm L $ we assume without loss of generality that $x_0=0$ and take $\delta$ sufficiently small, in such a way that, setting $x^\pm=\pm \delta$, we have: $d_h(x^\pm) \to 0$, $u_h(x^\pm)\to u(x^\pm) $. 

The points $x^\pm$ lie within mesh elements $e_h^\pm$ and for each element, we select a vertex $x_h^\pm$.  
Observe that, on each element $e_h$, as $h\to 0$, $d_h( x^r_h)-d_h(x^l_h)\to 0$, where $x^l_h$ and $x^r_h$
denote the left and right vertices of $e_h$, respectively; indeed, from the definition of the discrete energy $\F_h$, we have the estimate:
\begin{align}
    C&\geq \F_h( u_h, d_h,e_h) \geq\frac{G_c}{2}\frac{\epsilon_h}{h}(d_h(x^r_h)-d_h( x^l_h))^2
\end{align}
and since $h=o(\epsilon_h)$ as $h\to 0$, it follows that $d_h( x^r_h)-d_h(x^l_h)\to 0$. In particular, since $d_h(x^\pm)\to 0$, we get $d_h ( x^\pm_h) \to 0$ and $\bar d_h|_{e_h^\pm}\to 0$. For the boundary cases, if $x_0= L$ we set $x_h^-$ as above and $x_h^+=L$, while for $x_0=-L$, $x_h^-=-L$ and $x_h^+$ as above. Since the boundary conditions impose $ d_h(\pm L)= 0$, the same argument applies, ensuring that $d_h ( x^\pm_h) \to 0$ and $\bar d_h|_{e_h^\pm}\to 0$ also in this case.
  
We define $I_h=[x^-_h,x^+_h] $ and let $J_h=\{e_h \subset I_h\}$ denote the set of elements contained in $I_h$.
We then introduce  
 $\hat d_h=\max\{\bar d_h|_{e_h}\}$ and select an element $\hat e_h$ on which $ \bar d_h=\hat d_h$. Moreover, we call $J_h^\sharp=\{e_h \in J_h\,: \,
 u'_h|_{e_h}\geq a(\hat d_h)\frac{\sigma_c}{E_0}\} $ and accordingly we denote $I^\sharp_h$ the union of the elements $e_h \in J^\sharp_h$. Next, we define
 $$\bar u_h'=\begin{cases}
     u'_h & I_h\setminus ( I_h^\sharp\cup \hat e_h), 
     \\ a(\hat d_h)\frac{\sigma_c}{E_0} & I_h^\sharp\setminus  \hat e_h, 
      \\ \big( \sum_{I_h^\sharp\cup \,\hat{e}_h} u_h' \big) -\#(I^\sharp_h\setminus \hat e_h)\, a(\hat d_h)\frac{\sigma_c}{E_0} & \hat e_h. 
 \end{cases}$$
Note that $\bar u'_h\leq u'_h$ and $\bar{u}_h' \le a (\hat{d}_h) \frac{\sigma_c}{E_0}$ in $I_h\setminus\hat e_h$. Define $\bar u_h(x)=u_h(x_h^-)+\int_{x_h^-}^x\bar u_h'(r)\,dr$, and observe that
 $$u_h(x^+_h)-u_h(x_h^-)=\int_{I_h}\bar u'_h\text{d}x\leq \int_{\hat e_h}\bar u'_h\text{d}x+ \int_{I_h\setminus \hat e_h}a(\hat d_h)\frac{\sigma_c}{E_0}\text{d}x\leq h\bar u'_h+C\delta. $$
 As a consequence, in the element $\hat{e}_h$ for $h$ small enough and $\delta$ small enough we have 
\begin{equation} 
    \label{uBarra'}\bar u'_h\geq \frac{u_h(x^+_h)-u_h(x^-_h)-C\delta}{h} \ge a (\hat{d}_h) \frac{\sigma_c}{E_0} . 
\end{equation} 
Clearly,
\begin{equation}\label{emn1}
    \int_{I_h \setminus ( I_h^\sharp \cup \hat{e}_h)} f (u'_h, \bar{d}_h) \, \text{d}x 
    = \int_{I_h \setminus ( I_h^\sharp \cup \hat{e}_h)} f (\bar{u}'_h, \bar{d}_h) \, \text{d}x .
\end{equation} 
 Denoting $f(u_h,d_h,e_h)$ the restriction of $f(u_h,d_h)$ to the element $e_h$,  from \Cref{fBar}, it follows that 
\begin{align} \label{uBarra}
 \int_{I^\sharp_h \cup \hat{e}_h}f(u'_h,\bar d_h)\, \text{d}x & = h \sum_{e_h \in (J_h^\sharp \cup \hat{e}_h)}f(u'_h,\bar d_h, e_h) \\ & \geq h \Big( \sum_{e_h \in (J_h^\sharp \setminus \hat e_h)}f(\bar u'_h , \bar d_h, e_h) \Big) + h f(\bar u'_h, \hat d_h, \hat e_h) 
 \\ & \ge \int_{I^\sharp_h \setminus \hat{e}_h}  f (\bar u'_h , \bar d_h )  \, \text{d}x +  \int_{\hat{e}_h}  f(\bar u'_h, \hat d_h) \, \text{d}x .  \nonumber 
\end{align}
By \eqref{uBarra'}  $\bar u'_h$ exceeds the threshold $a(\hat{d}_h) \frac{\sigma_c}{E_0}$ on $\hat e_h$,  hence 
\begin{align}
    \int_{\hat e_h}f(\bar u'_h,\hat d_h)\, \text{d}x
    &\geq h a(\hat d_h)\sigma_c \bar{u}'_h - h \frac{\sigma_c^2}{2E_0}a^2(\hat d_h) 
    \\
    &= a(\hat d_h)\sigma_c j_h - o(1) , 
\end{align}
where $j_h=h \bar{u}'_h$.

Set $M_h=\max\{d_h(\hat x_h^r),d_h(\hat x_h^l)\}$ and $m_h=\min\{d_h(\hat x_h^r),d_h(\hat x_h^l)\}$, by the fact that $a$ is non-increasing and  Lipschitz continuous on $[0,1]$, it follows that $$a(\hat d_h)\geq a(M_h) \geq a(m_h)-2(M_h-m_h).$$
As a consequence, since for $h\to 0$, $d_h( \hat x^r_h)-d_h(\hat x^l_h)=o(1)$:
\begin{align} \label{adh}
    \int_{\hat e_h}f(\bar u'_h,\hat d_h)\text{d}x\geq  a(m_h)\sigma_cj_h-o(1)
\end{align}
We now focus on the remaining part of $\F_h(u_h,d_h,I_h)$ where, for $h$ small enough,
\begin{align} \label{zh}
  \tfrac{G_c}2\int_{I_h}\frac{1}{\epsilon_h}d_h^2+\epsilon_h|d_h'|^2\,\text{d}x&\geq\min\bigg\{\tfrac{G_c}2 \int_{I_h}\frac{1}{\epsilon_h}w_h^2+\epsilon_h|w_h'|^2\,\text{d}x \,:\,
w_h\in \mathbb{P}_{h}^1,\\&\qquad \qquad 
w_h(\hat x_h^l)=w_h(\hat x_h^r)=m_h,\,\, w_h(x_h^\pm)=d_h(x_h^\pm) \bigg\}\\&\geq\min\bigg\{
\tfrac{G_c}2 \int^{x_h^+/\epsilon_h}_{x_h^-/\epsilon_h}z_h^2+|z_h'|^2\,\text{d}x \,: \,
z_h\in \mathbb{P}_{\tilde h}^1%
,\\&\qquad \qquad z_h(\hat x_h^l/\epsilon_h)=z_h(\hat x_h^r/\epsilon_h)=m_h,\,\, z_h(x_h^\pm/\epsilon_h)=d_h(x_h^\pm) \bigg\},
\end{align}
where $\tilde{h} = h / \epsilon_h$ and by $\mathbb{P}_{\tilde h}^1$ we mean the piecewise affine functions on the mesh rescaled by $1/\epsilon_h$, that we call $\mathcal T_{\tilde{h}}$ and is defined on $\Big(\frac{x_h^-}{\epsilon_h},\,\frac{x_h^+}{\epsilon_h}\Big)$. 
We consider an extension of $z_h$ defined as follows:
$$  \tilde z_h(x)=\begin{cases}z_h(x_h^-/\epsilon_h)\Big(x-\frac{x_h^-}{\epsilon_h}+1\Big) & x\in\Big(\frac{x_h^-}{\epsilon_h}-1,\,\frac{x_h^-}{\epsilon_h}\Big)\\
z_h(x) & x\in \Big(\frac{x_h^-}{\epsilon_h},\,\frac{x_h^+}{\epsilon_h}\Big)
   \\z_h(x_h^+/\epsilon_h)\Big(\frac{x_h^+}{\epsilon_h}+1-x\Big) & x\in\Big(\frac{x_h^+}{\epsilon_h},\,\frac{x_h^+}{\epsilon_h}+1\Big)
\end{cases}$$
and observe that $$\int_{x_h^+/\epsilon_h}^{x_h^+/\epsilon_h+1}\tilde z_h^2+|\tilde z_h'|^2\,\text{d}x=z_h^2(x_h^+/\epsilon_h)\int_{[0,1]}(1-x)^2+1\,\text{d}x=c z_h^2(x_h^+/\epsilon_h) = c d_h^2(x^+_h) \to 0.$$
The same reasoning can be applied to $\Big(\frac{x_h^-}{\epsilon_h}-1,\,\frac{x_h^-}{\epsilon_h}\Big)$ and hence 
\begin{align} \label{ch}
  \int_{x_h^-/\epsilon_h}^{x_h^+/\epsilon_h} z_h^2+| z_h'|^2\,\text{d}x&= \int_{x_h^-/\epsilon_h-1}^{x_h^+/\epsilon_h+1}\tilde z_h^2+|\tilde z_h'|^2\,\text{d}x-o(1).
\end{align}
For $h$ small enough, we can therefore focus our analysis on the study of the optimal profile of the functions $\tilde z_h\in \mathbb{P}_{\tilde h}^1$\, on the interval $\Big(\frac{x_h^-}{\epsilon_h}-1,\,\frac{x_h^+}{\epsilon_h}+1\Big)$ such that   $ \tilde z_h|_{\hat e_h}=m_h,$ and $\tilde z_h(x_h^\pm/\epsilon_h\pm 1)=0$. 
We introduce the localized energies
$K_R(z)=\int_{(0,R)}z^2+|z'|^2\text{d}x$ and call $z_{R,h}$ the solutions of the minimization problems:
$$z_{R,h}\in\argmin\{K_R(z)\,:\,z\in H^1(0,R),\,\,z(0)=m_h,\,z(R)=0\}.$$
We call $R_h^+=\frac{x_h^+}{\epsilon_h}+ 1-\frac{\hat x_h^r}{\epsilon_h}$ and $R_h^-=\frac{\hat x_h^l}{\epsilon_h}-\Big(\frac{x_h^-}{\epsilon_h}- 1\Big)$ and observe that, set $\tilde R_h=\max\{R_h^\pm\}$, 
\begin{equation} \label{Rtilde}
   K_{\tilde R_h}(z_{\tilde R_h,h})\leq K_{R_h^\pm}(z_{R_h^\pm,h}). 
\end{equation} 
Therefore
\begin{align} 
 \int_{x_h^-/\epsilon_h-1}^{x_h^+/\epsilon_h+1} \tilde z_h^2+|\tilde z_h'|^2\,\text{d}x\geq K_{R_h^-}(z_{R_h^-,h})+ K_{R_h^+}(z_{R_h^+,h})\geq  2 K_{\tilde R_h}(z_{\tilde R_h,h})\geq 2m_h^2,
\end{align} 
where the last inequality follows from the study of the optimal profile in \Cref{optprof}.
Combining this estimate to \eqref{zh} and \eqref{ch} leads to:
\begin{align} \label{fracinf}
\tfrac{G_c}2\int_{I_h}\frac{1}{\epsilon_h}d_h^2+\epsilon_h|d_h'|^2\,\text{d}x
\geq G_c m_h^2- o(1).
\end{align} 
Taking the sum \eqref{fracinf} and \eqref{adh} we obtain 
\begin{align} 
\int_{\hat{e}_h} f(\bar{u}_h',\bar d_h)\text{d}x + \tfrac{G_c}2\int_{I_h}\frac{1}{\epsilon_h}d_h^2+\epsilon_h|d_h'|^2\,\text{d}x 
&\geq a(m_h)\sigma_cj_h + G_c m_h^2 -o(1) 
\\&\geq \phi(j_h)-o(1) .
\end{align}
Hence, 
\begin{align} 
\F_h(u_h,d_h,I_h)&=\int_{I_h} f(u_h',\bar d_h)\text{d}x+\tfrac{G_c}2\int_{I_h}\frac{1}{\epsilon_h}d_h^2+\epsilon_h|d_h'|^2\,\text{d}x 
\\&\geq \int_{I_h\setminus \hat{e}_h} f(\bar{u}_h',\bar d_h)\text{d}x + \phi (j_h)-o(1).
\end{align}
Let $\hat{e}_h = (\hat{x}^-_h, \hat{x}^+_h)$ and define locally in $I_h$ the function
\begin{equation}
\label{hatuh}\hat u_h(x)=\begin{cases}
    \bar{u}_h(x) & x \in I_h \setminus \hat{e}_h , 
    \\ \bar{u}_h(\hat{x}_h^-)\mathbbm{1}_{[\hat x_h^-, \hat x_h)}(x) 
     +u( \hat{x}_h^+)\mathbbm{1}_{[\hat x_h, \hat x_h^+ ]}(x)  & x\in \hat{e}_h ,
     \end{cases}    
\end{equation}
     where $\hat{x}_h \in (\hat{x}_h^- , \hat{x}_h^+)$. 
Then, by \eqref{uBarra'} $ \llbracket \hat{u}_h \rrbracket = h \bar{u}_h' = j_h$
and thus 
$$
   \F_h(u_h,d_h,I_h) \ge \int_{I_h} f( \hat{u}_h', \bar d_h)\text{d}x + \phi_n( \llbracket \hat{u}_h \rrbracket)-o(1) .
$$
Note that in $I_h \setminus \hat{e}_h$ we have $\hat{u}_h' = \bar{u}_h' \le a (\hat{d}_h) \frac{\sigma_c}{E_0}$ while $\hat{u}_h'=0$ in $\hat{e}_h$. Hence, by \Cref{l.fLip} 
$$
    f ( \hat{u}_h' , \bar{d}_h ) \ge f(\hat{u}_h', 0) - C | \bar{d}_h| | \hat{u}'_h| \ge f ( \hat{u}_h' , 0) - C' \ge f_n ( \hat{u}_h' , 0) - C' . 
$$
In conclusion
\begin{equation} \label{mn1}
   \F_h(u_h,d_h,I_h) \ge \int_{I_h} f_n ( \hat{u}_h', 0 ) \, \text{d}x + \phi_n( \llbracket \hat{u}_h \rrbracket)-o(1) - C \delta.
\end{equation}

We apply the same reasoning to each point $x_i\in J_u^t$ for $i=0,...,N_t = \# J_u^t$. 
We define $x_{i,h}^\pm$ analogously to $x_h^\pm$ and then we set $I^i_h=[x_{i,h}^-,x_{i,h}^+]$ and 
$J^{t,\delta}_{u,h}=\bigcup_{i=1}^{N_t}I_h^i$. Defining $\hat{u}_h$ in each interval $I^i_h$ as above  and summing over $i=1,...,N_t$, we obtain:
\begin{align} \label{infphi}\F_h\big(u_h,d_h, J^{t,\delta}_{u,h}\big)\geq \int_{J_{u,h}^{t,\delta}}   f_n ( \hat{u}_h' , 0 ) \, \text{d}x + \sum_{i=1}^{N_t} \phi_n( \llbracket \hat{u}_h (\hat{x}_{i,h}) \rrbracket ) - o(1) - C \delta . 
\end{align}

\medskip
{\bf Out of jump points.} 
Since $d_h\to 0$ in $L^1(I)$, it also converges quasi-uniformly, namely for $\epsilon>0$, there exists $I_\epsilon\subset I$ such that $|I_\epsilon|<\epsilon$ and $d_h \to 0$ uniformly on $I\setminus I_\epsilon$.
If we restrict to $I\setminus (J_{u,h}^{t,\delta} \cup I_\epsilon)$, where we have uniform convergence, for every $\gamma >0$, there exists $h_\gamma$ such that for every $h<h_\gamma$, $\|d_h\|_{L^\infty(I\setminus I_\epsilon)}\leq \gamma$ and therefore:
 \begin{equation}
     \F_h(u_h,d_h,I\setminus (J_{u,h}^{t,\delta}\cup I_\epsilon))\geq \int_{I\setminus (J_{u,h}^{t,\delta} \cup I_\epsilon)}f(u'_h,\gamma) \,\text{d}x
 \end{equation}
By convexity, $f(s,r) \ge a(r) \sigma_c |s| - \frac{\sigma_c^2}{2E_0}a^2(r)$. It follows that 
$$\int_{I\setminus (J_{u,h}^{t,\delta} \cup I_\epsilon)}f(u'_h,\gamma) \,\text{d}x \ge 
\int_{I\setminus (J_{u,h}^{t,\delta} \cup I_\epsilon)} a(\gamma) \sigma_c | u'_h| \, \text{d}x - C |I\setminus (J_{u,h}^{t,\delta} \cup I_\epsilon)| . $$
Hence $u'_h$ is bounded in $L^1(I \setminus (J_{u,h}^{t,\delta} \cup I_\epsilon))$. By Lemma \ref{l.fLip} we have $f (s,r) \ge f (s,0) - C | r | | s|$. Thus
\begin{align} \label{inff}
\F_h(u_h,d_h,I\setminus (J_{u,h}^{t,\delta}\cup I_\epsilon)) & \geq
\int_{I\setminus (J_{u,h}^{t,\delta} \cup I_\epsilon)}f (u'_h,\gamma)\,\text{d}x \\ & \ge \int_{I\setminus (J_{u,h}^{t,\delta} \cup I_\epsilon)}f (u'_h,0)\,\text{d}x - C | \gamma| \int_{I\setminus (J_{u,h}^{t,\delta} \cup I_\epsilon)} |u'_h|\, \text{d}x  \\
& \ge \int_{I\setminus (J_{u,h}^{t,\delta} \cup I_\epsilon)}f_n(u'_h,0)\,\text{d}x - C' \gamma. 
\end{align}
We define $\hat{u}_h = u_h$ in $I \setminus J^{t,\delta}_{u,h}$. Then, 
taking the sum of \eqref{infphi} and \eqref{inff} and recalling \eqref{hatuh}, 
we obtain:
\begin{align} \label{liminf1}
\F_h (u_h,d_h) &
\geq \F_h (u_h,d_h, I\setminus I_\epsilon ) \\
& \geq\int_{I\setminus I_\epsilon }f_n(\hat u'_h,0)\text{d}x +\sum_{x_i \in J_{\hat{u}}} \phi_n(\llbracket \hat u_h (x_i)\rrbracket  )-o(1)- C \delta - C' \gamma \\&=\G_n (\hat u_h,I\setminus I_\epsilon) -o(1)- C \delta - C' \gamma .
 \end{align}

\medskip
 {\bf Liminf.} We take the liminf of both sides of the previous inequality. Since $u_h \rightharpoonup u$ in $\BV(I)$, it follows that (up to subsequences) $\hat u_h$ converges weakly in $\BV(I)$ to a certain function $u_{t,\delta}$. By the definition of $\hat{u}_h$ we obtain $u_{t,\delta}=u$ in $I \setminus J_u^{t,\delta}$, where $J_u^{t,\delta}$ is the union of intervals $I_i = [ x_i - \delta , x_i + \delta]$ for $x_i \in J_u^t$. In each interval $I_i$, the function $u_{t,\delta}$ has a jump, in a certain point $\hat{x}_i$, with 
 $$u_{t,\delta}^+ (\hat{x}_i) - u_{t,\delta}^- (\hat{x}_i) \ge u(x_i^+) - u(x_i^-) - C \delta $$ since by \eqref{uBarra'}
 $$
     \hat{u}_h (\hat{x}_h^+) - \hat{u}_h (\hat{x}_h^-) = h \bar{u}'_h (\hat{x}_h) \ge u_h (x_h^+) - u_h (x_h^-) - C \delta \ \to \ u (x_i^+) - u (x_i^-) - C \delta . 
 $$
 We recall that $x_i^\pm$ denotes $x_i\pm \delta$ and therefore $u (x_i^+) - u(x_i^-) \to u^+ (x_i) - u^- (x_i)$ as $\delta \to 0$.
 
Therefore, recalling Corollary \ref{relaxGn}, we obtain:
 \begin{align}
 \liminf_{h\to 0} \F_h (u_h,d_h)
 & \geq \liminf_{h\to 0}\G_n (\hat u_h, I\setminus I_\epsilon) - C \delta - C' \gamma \\
& \geq \bar\G_n (u_{t,\delta}, I\setminus I_\epsilon) - C \delta - C' \gamma 
 \end{align}
 for every $n\in\mathbb{N}$, $t,\delta,\epsilon,\gamma>0$. Taking the supremum with respect to $\epsilon$ and $\gamma$ yields 
$$
    \liminf_{h\to 0} \F_h (u_h , d_h) \ge \bar\G_n (u_{t,\delta}) - C \delta . 
$$
Taking the supremum with respect to $n$ and recalling \Cref{supGn}, yields:
 \begin{align}
 \liminf_{h\to 0}\F_h
(u_h,d_h)\geq \sup_{n\in\mathbb{N}}\bar\G_n(u_{t,\delta})-   C \delta  =\F(u_{t,\delta})  - C \delta  . 
 \end{align}
It remains to pass to the limit with respect to $t$ and $\delta$. To this end, being $u_{t,\delta}=u$ in $I \setminus J_u^{t,\delta}$ we can write 
\begin{align}
    \F ( u_{t,\delta} ) & = \F ( u_{t,\delta} , I \setminus J_u^{t,\delta} ) + \F ( u_{t,\delta} ,J _u^{t,\delta} )
    \\ & \ge \F ( u , I \setminus J_u^{t,\delta}  ) + \sum_{x \in J_u^t} \phi (\llbracket u_{t,\delta}(x) \rrbracket) . 
\end{align} 
As stated above 
 $$ \llbracket u_{t,\delta} (\hat{x}_i) \rrbracket \ge u(x_i^+) - u(x_i^-) - C \delta \ \to \ \llbracket u (x_i) \rrbracket \quad \text{as $\delta \to 0$.}  $$ 
Moreover $I \setminus J_u^{t,\delta} \nearrow I \setminus J_u^t$ as $\delta \to 0$. Hence
$$
    \liminf_{h\to 0} \F_h (u_h,d_h) \geq \F ( u , I \setminus J_u^{t}  ) +  \sum_{x \in J_u^t} \phi (\llbracket u (x)  \rrbracket) . 
$$
Taking the supremum with respect to $t>0$ yields
$$
    \liminf_{h\to 0} \F_h (u_h,d_h) \geq \F ( u , I \setminus J_u ) +  \sum_{x \in J_u} \phi (\llbracket u (x)  \rrbracket) ,
$$
which concludes the proof. \qed

\section{Numerical results\label{s.num}}

This section presents numerical results that validate the isotropy of the discrete energy.  Specifically, we assess whether the numerical approximation introduces any mesh-induced anisotropy and show that the formulation remains robust with respect to the mesh geometry. These results confirm that displacement jumps do not compromise the isotropic character of the fracture energy, which is governed by the phase-field profile over a neighborhood of size $ \epsilon_h \gg h $. 
To this end, we consider the multiaxial energy employed in \cite{vicent} under plane-strain conditions. Specifically, the energy takes the form
$$
    F_h ( \bm{u}_h , \bm{\eta}_h , d_h ) = \int_\Omega \psi_e ( \bm{\varepsilon}_h - \bm{\eta}_h   ) \, \text{d}x + \int_\Omega \pi (\bm{\eta}_h, d_h) \, \text{d}x + \tfrac{G_c}{2} \int_\Omega \frac{1}{\eps_h} d_h^2 + \eps_h | \nabla d_h |^2 \, \text{d}x .
$$

In general, as detailed in \cite{vicent}, the onset of material nonlinearities is governed by the specification of the elastic domain, within which the stress tensor $\bm{\sigma}$ is constrained to lie. 
When the eigen-strain $\bm{\eta}_h$ becomes non-zero, nonlinear dissipative effects emerge. As further shown in \cite{vicent}, the eigen-strain potential $\pi(\cdot,
d_h)$ coincides with the support function of the elastic domain for a fixed value of the damage variable $d_h$.

Here, using the standard volumetric-deviatoric decomposition, the elastic energy density $\psi_e$ depends on the traces and deviatoric norms of the elastic strain tensor $(\bm{\strain}_h - \bm{\eta_h}$) and reads 
\begin{align}
\psi_e(\bm{\varepsilon}_h - \bm{\eta}_h) &= \hat{\psi}_e\big(\operatorname{tr}(\bm{\varepsilon}_h - \bm{\eta}_h), \|\bm{\varepsilon}_{h,\text{dev}}-\bm{\eta}_{h,\text{dev}}\|\big) \\
&:= \frac{\kappa}{2} \operatorname{tr}^2(\bm{\varepsilon}_h-\bm{\eta}_h)  
+ \mu  \|\bm{\varepsilon}_{h,\text{dev}}-\bm{\eta}_{h,\text{dev}}\| ^2,
\end{align}
where $\kappa$ and $\mu$ are the bulk and shear moduli, respectively. 

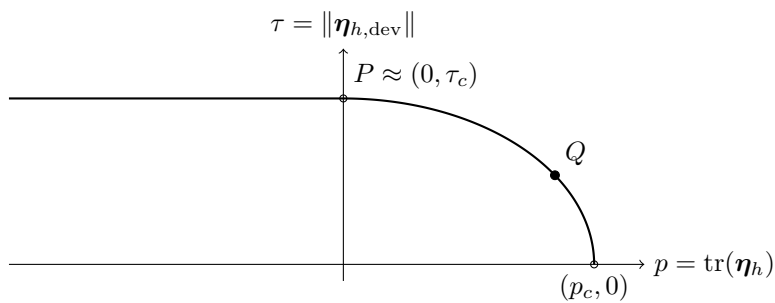
\begin{figure}[ht] 
    \centering
    \begin{tikzpicture}[scale=1.1, every node/.style={font=\small}]
      \def\pc{3}   
      \def\tc{2}   
      \def\ang{23} 

      \pgfmathsetmacro{\m}{tan(\ang)}

      \pgfmathsetmacro{\px}{1/sqrt(1/(\pc*\pc) + (\m*\m)/(\tc*\tc))}
      \pgfmathsetmacro{\py}{\m*\px}

      \draw[->] (-4,0) -- ({\pc+0.6},0) node[right] {$p=\operatorname{tr}(\bm\eta_h)$};
      \draw[->] (0,-0.2) -- (0,{\tc+0.6}) node[above] {$\tau=\|\bm\eta_{h,\mathrm{dev}}\|$};

      \draw[thick, black, samples=200, smooth]
        plot[domain=0:90] ({\pc*cos(\x)}, {\tc*sin(\x)});
      \draw[thick, black] (-4,\tc) -- (0,\tc);

      \draw (\pc,0) circle (1.2pt) node[below] {$(p_c,0)$};
      \draw (0,\tc) circle (1.2pt) node[above right] {$P\approx(0,\tau_c)$};

      \filldraw (\px,\py) circle (1.5pt) node[above right] {$Q$};
    \end{tikzpicture}
    \caption{Shape of the elastic domain. Points P and Q are the stress states at which the material starts fracturing for the two tests discussed below.}
    \label{domain}
\end{figure}

As discussed above and detailed in \cite{vicent}, the model is able to reproduce a variety of strength surfaces consistent with experimental data. In the following, we present the numerical results obtained using the eigenstrain potential
\[
\pi(\bm{\eta}_h,d_h) = 
\begin{cases}
a(d_h)\cdot\phi_2(\operatorname{tr}(\bm{\eta}_h) , \|\bm{\eta}_{h,\text{dev}}\|), & \text{if } \operatorname{tr}(\bm{\eta}_h) \geq 0 \\
+\infty & \text{otherwise},
\end{cases}
\]
where
\[
\phi_2(\operatorname{tr}(\bm{\eta}_h), \|\bm{\eta}_{h,\text{dev}}\|) = \sqrt{p_c^2 \operatorname{tr}^2(\bm{\eta}_h) + \tau_c^2 \|\bm{\eta}_{h,\text{dev}}\|^2}.
\]
Such a potential defines a semi-elliptic strength surface that passes through $(p_c,0)$ and $(0,\tau_c)$, where $p_c$ and $\tau_c$ are respectively the critical pressure and shear stress. 
The shape of the elastic domain is depicted in Figure \ref{domain}.

\medskip

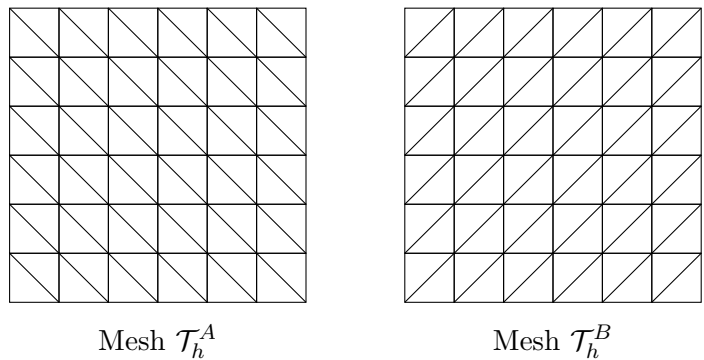
\begin{figure}[H] 
    \centering
    \begin{tikzpicture}[scale=0.65]  

        \def\N{6}        

        \foreach \i in {0,...,5} {
            \foreach \j in {0,...,5} {
                \draw (\i,\j) -- ++(1,0) -- ++(0,1) -- ++(-1,0) -- cycle;
                \draw (\i,\j+1) -- (\i+1,\j);
            }
        }

        \foreach \i in {0,...,5} {
            \foreach \j in {0,...,5} {
                \pgfmathtruncatemacro{\x}{\i + 8}
                \draw (\x,\j) -- ++(1,0) -- ++(0,1) -- ++(-1,0) -- cycle;
                \draw (\x,\j) -- (\x+1,\j+1);
            }
        }

        \node at (3,-0.8) {Mesh $\mathcal{T}_h^A$};
        \node at (11,-0.8) {Mesh $\mathcal{T}_h^B$};

    \end{tikzpicture}
    \caption{The two different mesh discretizations.}
    \label{meshQ}
\end{figure}

Simulations were performed using the \texttt{GRIPHFiTH} Matlab library for phase-field fracture modeling, available at \texttt{https://gitlab.ethz.ch/compmech}. 
We analyze an initially intact square domain with edge length $L=1$. For the material parameters, we set  $E_0= 10^3$, 
$\nu = 0.3$, $G_c = 0.2$, $\epsilon_h = 0.025$, and  $p_c = \tau_c = 10$. 
The domain is discretized using two distinct triangular meshes, each composed of right-angled triangles with leg length $h$, such that  $\epsilon_h/h  \approx 5$ (see Figure \ref{meshQ}).  

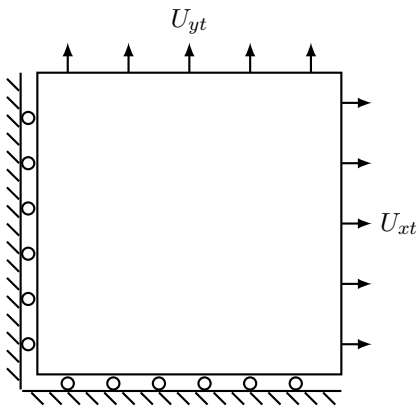
\begin{figure}[h]
\centering
\begin{tikzpicture}[scale=4]

\draw[thick] (0,0) rectangle (1,1);

\foreach \x in {0.1,0.3,0.5,0.7,0.9} {
    \draw[-latex, thick] (\x,1) -- ++(0,0.1);
}

\foreach \y in {0.1,0.3,0.5,0.7,0.9} {
    \draw[-latex, thick] (1,\y) -- ++(0.1,0);
}

\draw[thick] (-0.05,-0.055) -- (1,-0.055); %
\draw[thick] (-0.055,-0.05) -- (-0.055,1); 

\foreach \x in {0.1,0.25,0.4,0.55,0.7,0.85} {
    \draw[thick] (\x,-0.03) circle (0.02);
}

\foreach \y in {0.1,0.25,0.4,0.55,0.7,0.85} {
    \draw[thick] (-0.03,\y) circle (0.02);
}

\foreach \y in {0.02,0.08,...,0.98} {
    \draw[thick] (-0.1,\y) -- (-0.06,\y-0.04);
}

\foreach \x in {0.02,0.08,...,0.98} {
    \draw[thick] (\x,-0.1) -- (\x-0.04,-0.06);
} 

\node[scale=0.9] at (0.5,1.17) { $U_{yt}$};
\node[scale=0.9
] at (1.19,0.5) {$U_{xt}$};

\end{tikzpicture}
\caption{Set-up of the numerical simulations.}
\label{fig:setup2D}
\end{figure}

The boundary conditions are defined by enforcing $d_h=0$
along the entire boundary. Roller supports are applied along the left and bottom edges, allowing displacements only in the tangential direction. Normal displacements $U_{xt}$ and $U_{yt}$ are imposed on the right and top edges 
and are increased linearly over 1000 
 loading steps. 
The values of the imposed displacement at the final time step are called $U_x$ and $U_y$ respectively. The setup is illustrated in Figure \ref{fig:setup2D} .

Initially, the strain field remains homogeneous, with components $\epsilon_{xx}=U_{xt},\,\,\epsilon_{yy}=U_{yt}$ and $\epsilon_{xy}=0$. 
By varying the ratio between the imposed displacements on the top and right edges, the full range of stress states at which the material fractures can be explored (see Figure \ref{domain}). 

\begin{figure}
\centering
\includegraphics[width=.4\textwidth]{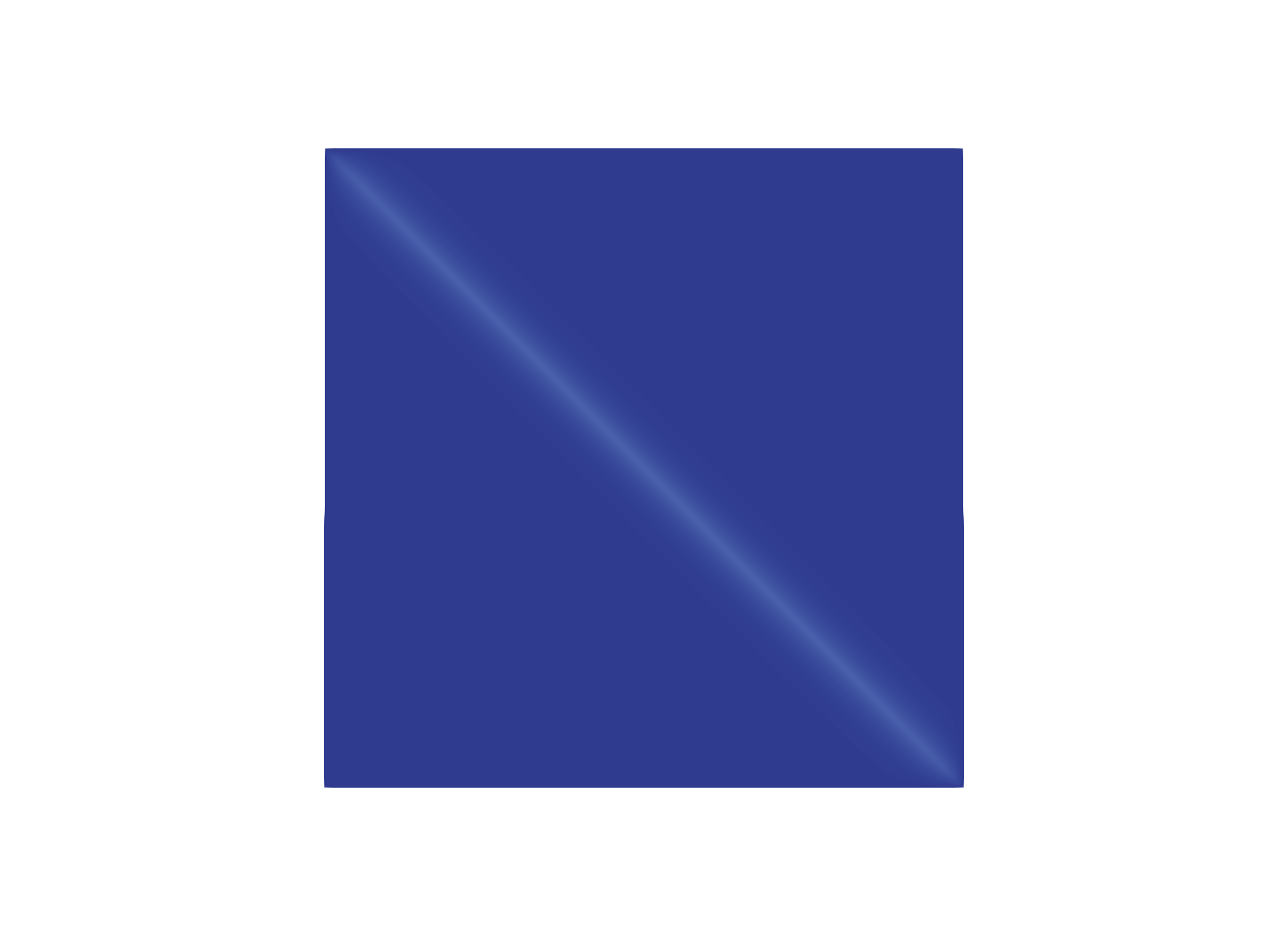}
\includegraphics[width=.4\textwidth]{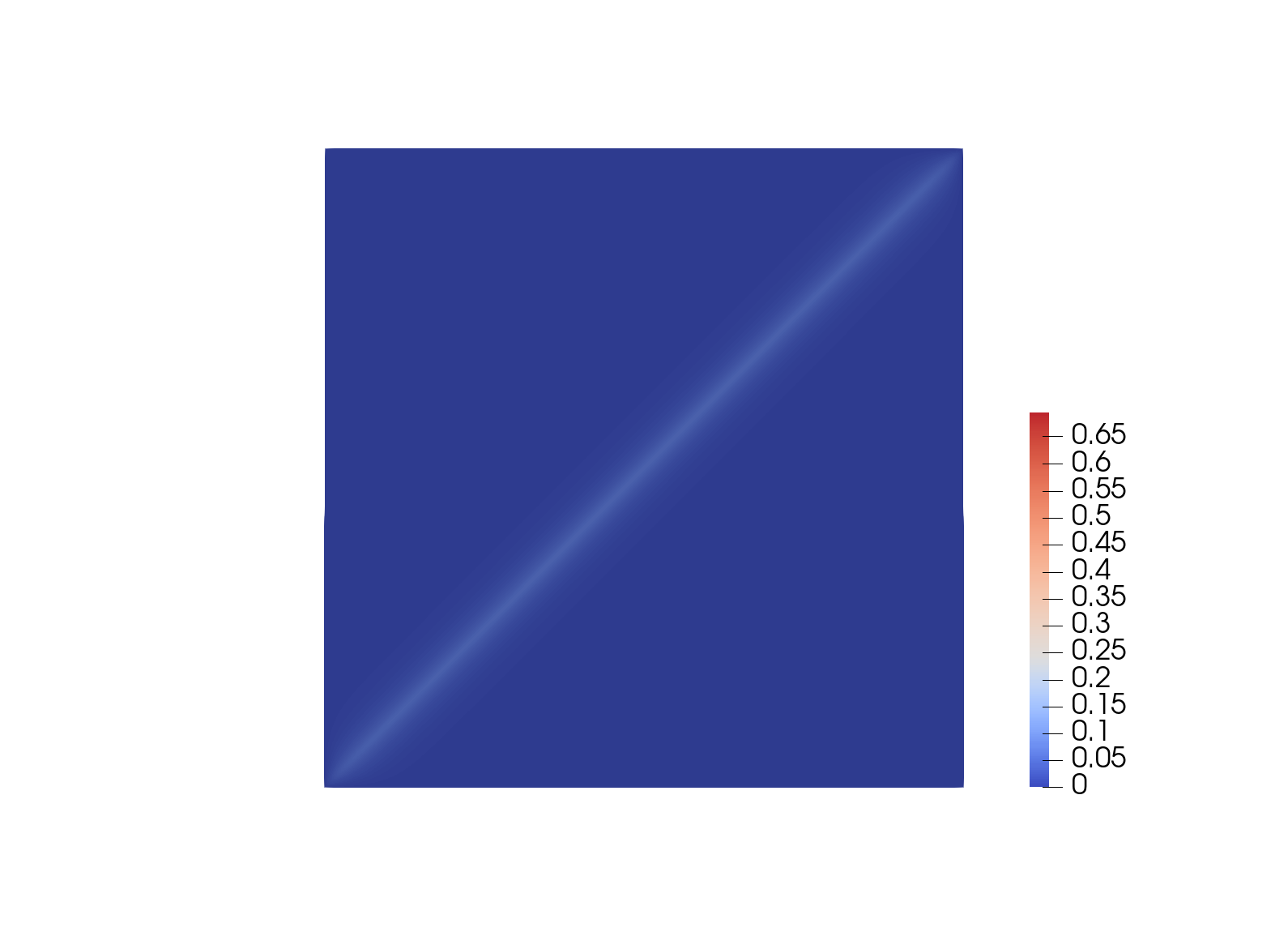}
\caption{
$(U_x,U_y)=(0.5,-0.45)$: 
loading configuration with multiple solutions. Phase fields for mesh $\mathcal T_A$ (left) and mesh $\mathcal T_B$ (left) at $(U_x,U_y)= (0.01,0.009)$ 
. In both cases it is localized on a strip of width $\epsilon_h=0.025.$}\label{fig:shearPF}
\end{figure}

Under certain loading configuration, the problem admits multiple solutions. For example, in the case of pure shear, that can be obtained setting $U_x\sim-U_y$, either diagonal may serve as a failure path. Damage localizes instantaneously along one of these patterns once the stress state reaches point P in Figure \ref{domain}. In this regime, the mesh topology influences which of these admissible solutions is selected by the algorithm, highlighting the sensitivity of non-unique solutions to the mesh (see Figure \ref{fig:shearPF}). 
Despite the different crack patterns, the damage localizes in both cases 
when $(U_x,U_y)= (0.01,0.009)$ \bl and the fracture energy at that time is the same and equal to $0.8309\cdot 10^{-3}$. 

\begin{figure}
\centering
\includegraphics[width=.4\textwidth]{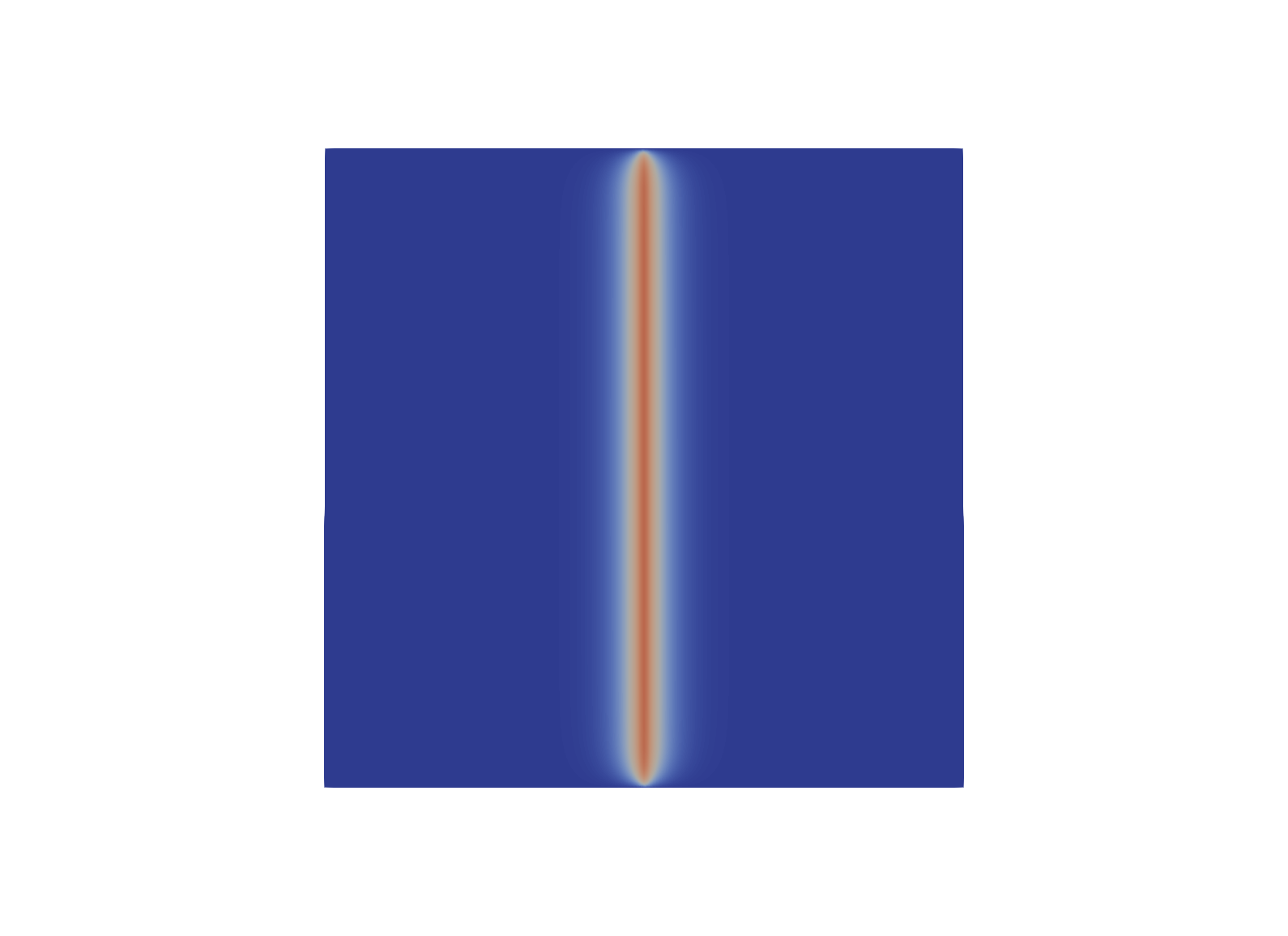}
\includegraphics[width=.4\textwidth]{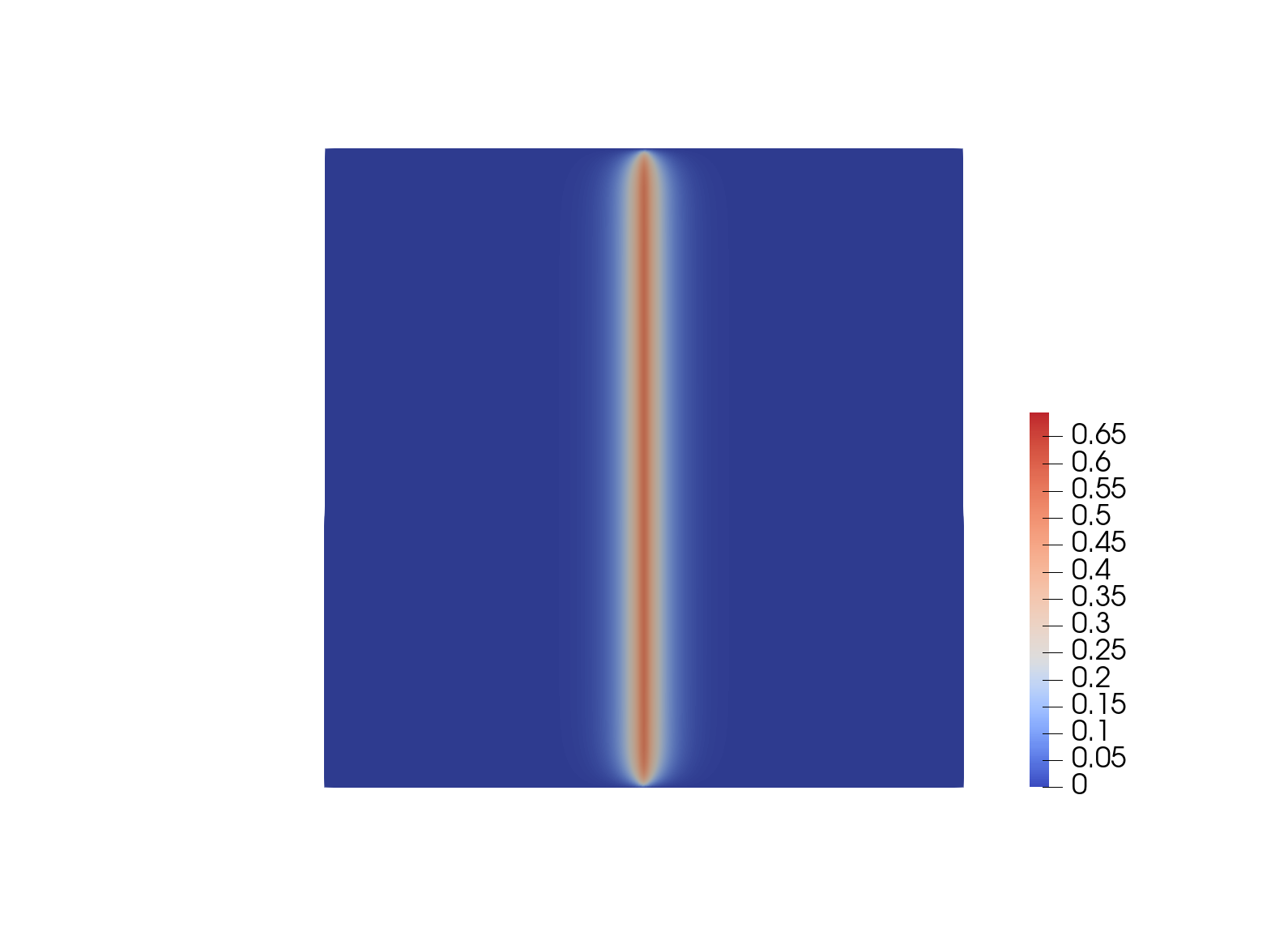}

\caption{
$(U_x,U_y)=(1,0.5)$: loading configuration with unique solution. Phase fields for mesh $\mathcal T_A$ (left) and mesh $\mathcal T_B$ (left) at $(U_x,U_y)=(0.017,0.0085)$. In both cases it is localized on a strip of width $\epsilon_h=0.025.$}\label{fig:vert}
\end{figure}

Naturally, if instead we perform a test 
that has as unique solution a vertical crack at the middle of the domain, the mesh has no influence on the result. This can be seen in Figure \ref{fig:vert}.   Such a result is obtained by prescribing the displacements $(U_x,U_y)=(1,0.5)$,  which corresponds to a failure stress state lying on the elastic domain at an angle $\theta \approx 23^\circ $  (point Q in Figure \ref{domain}).  The damage localizes in both cases when $(U_x,U_y)= (0.017,0.0085)$ and the fracture energy at that time is the same and equal to $0.0579106$. 
    
 \begin{remark}
In the proof of $\Gamma$-convergence, we exploited the tendency of the strain to concentrate in narrow regions. This behavior is confirmed numerically 
by Figure \ref{fig:strains}, where it is shown that the strain localizes within a narrow band of thickness proportional to the mesh size. On the other hand, the damage variable is distributed over a wider region proportional to the internal length $\epsilon_h$, as observable in Figures \ref {fig:shearPF} and \ref{fig:vert}. This numerical observation is consistent with the theoretical framework discussed in the previous sections: 
the crack is geometrically approximated by a narrow band of elements exhibiting large displacement gradients, while the regularization induced by $\epsilon_h$ ensures a smooth, mesh-independent damage profile. As a result, the fracture energy is captured through non-local contributions, preventing mesh bias. 
\begin{figure}
\centering
\includegraphics[width=.4\textwidth]{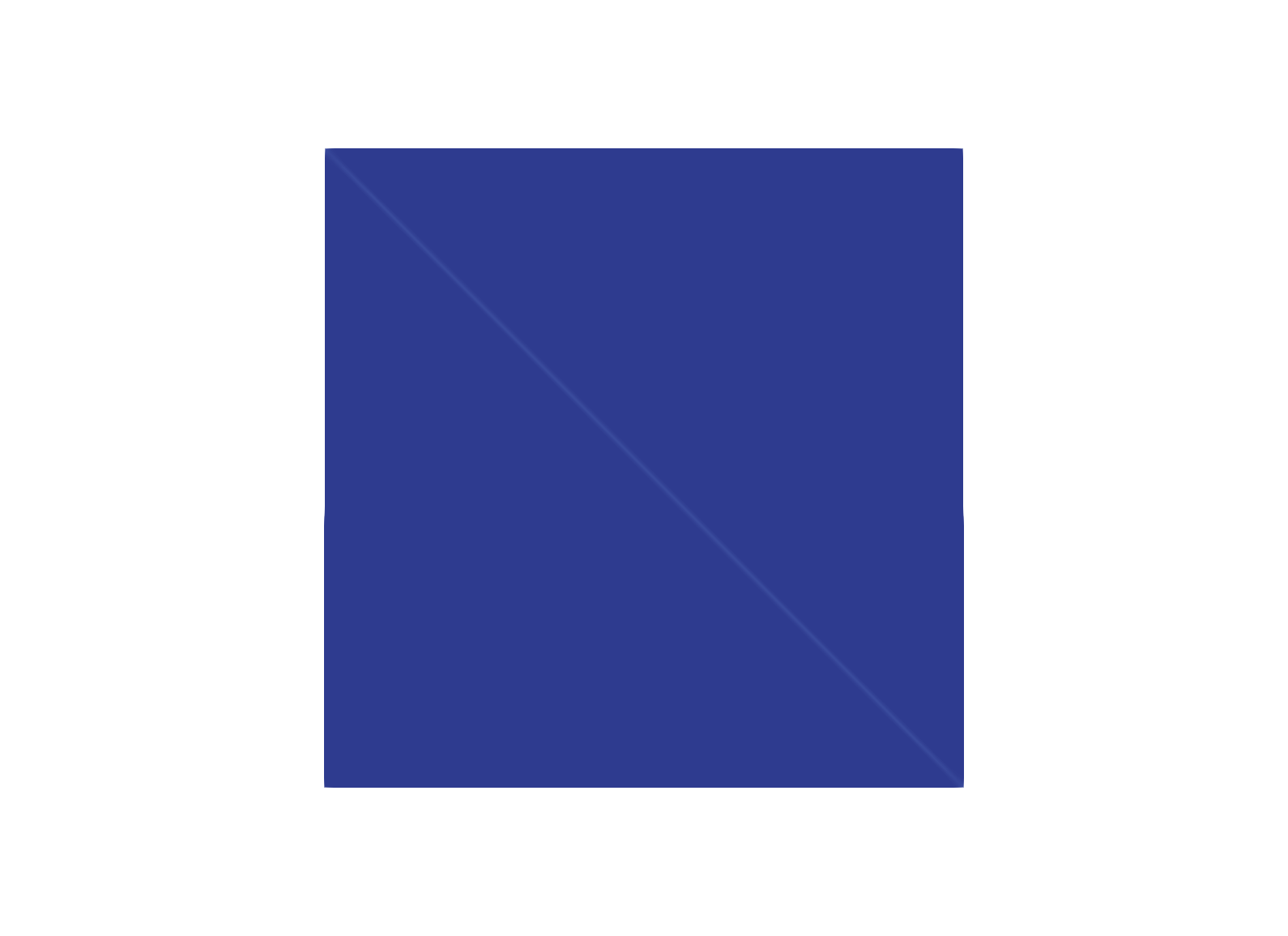}
\includegraphics[width=.4\textwidth]{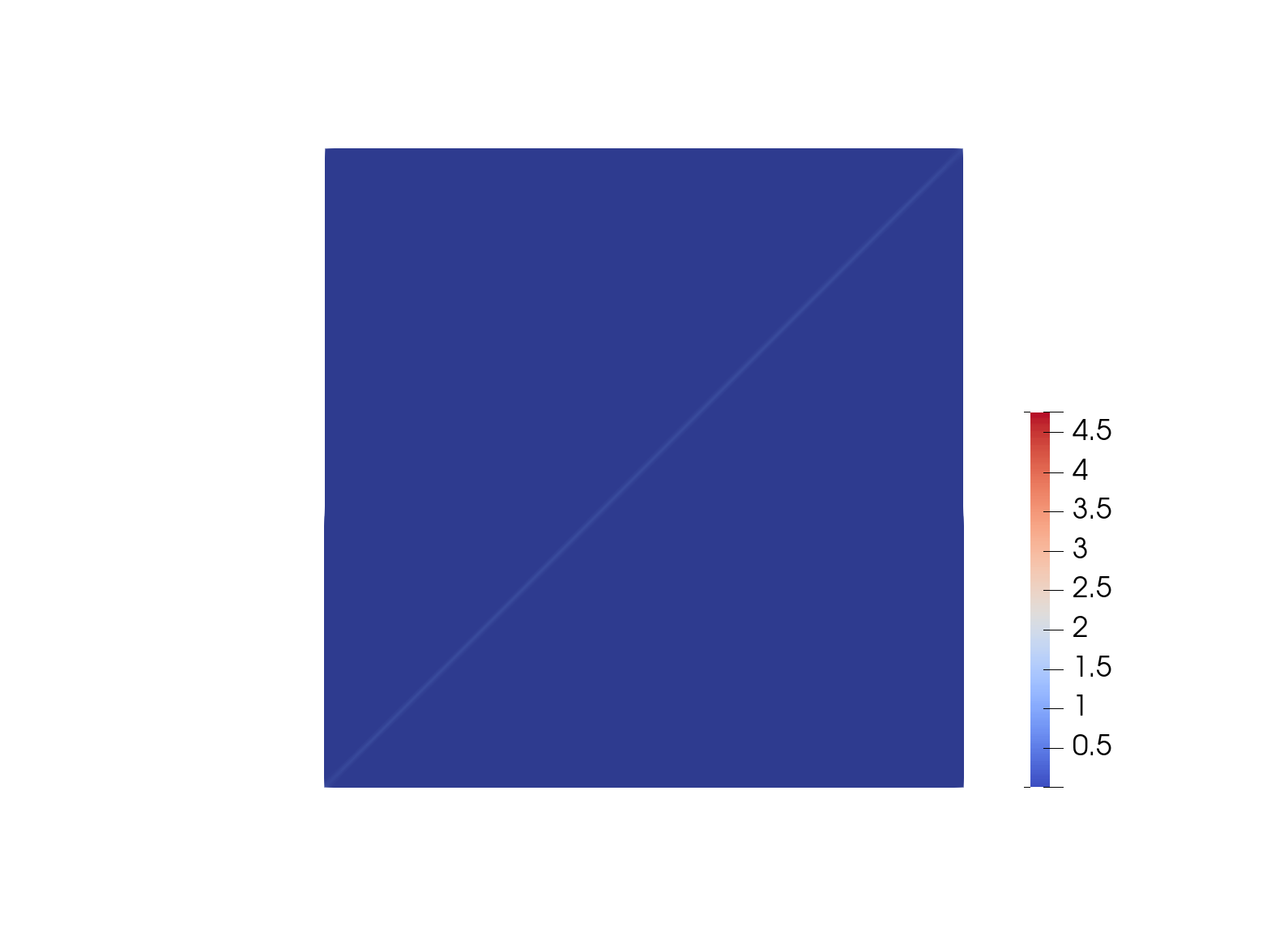}

\includegraphics[width=.4\textwidth]{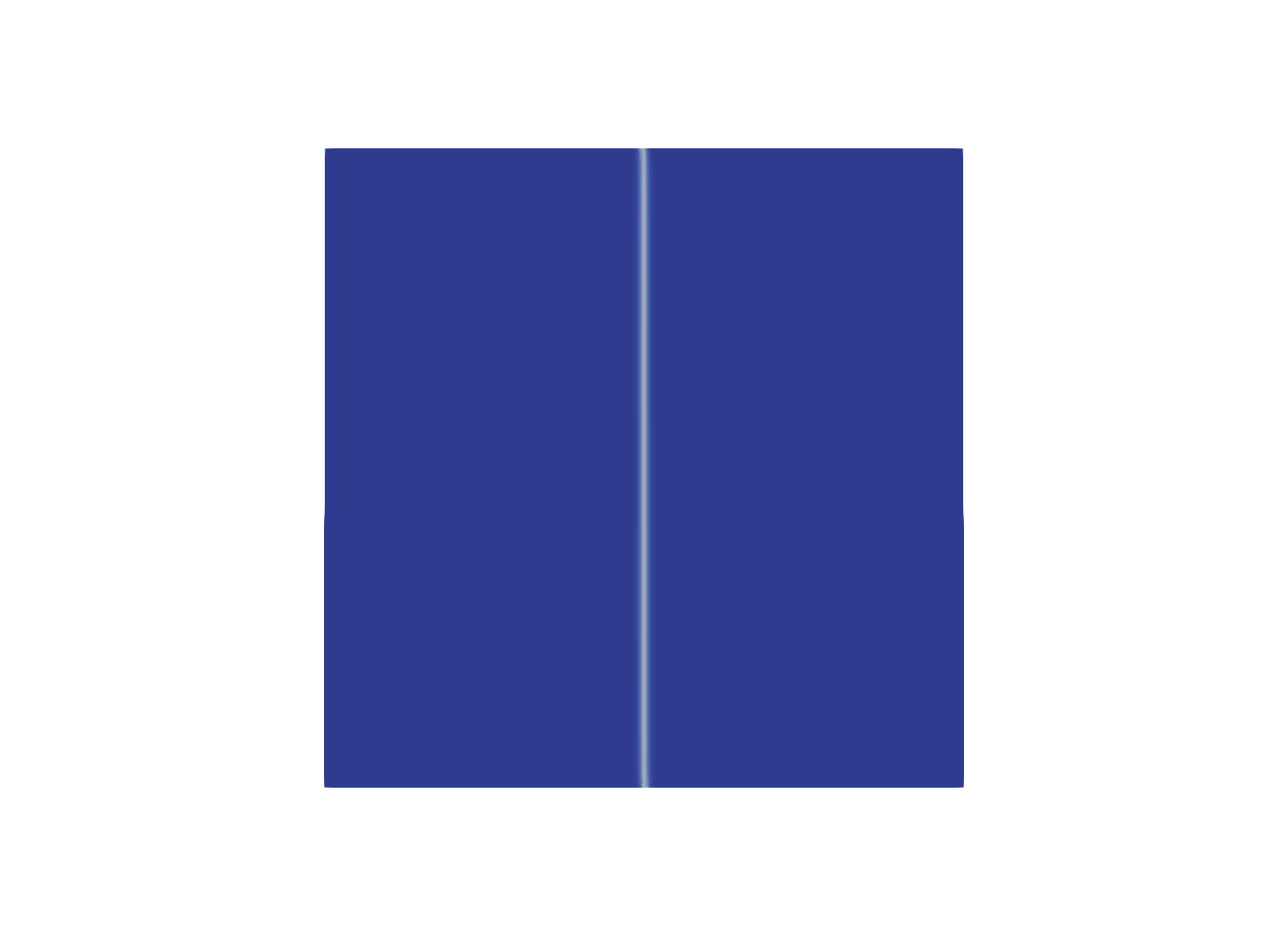}
\includegraphics[width=.4\textwidth]{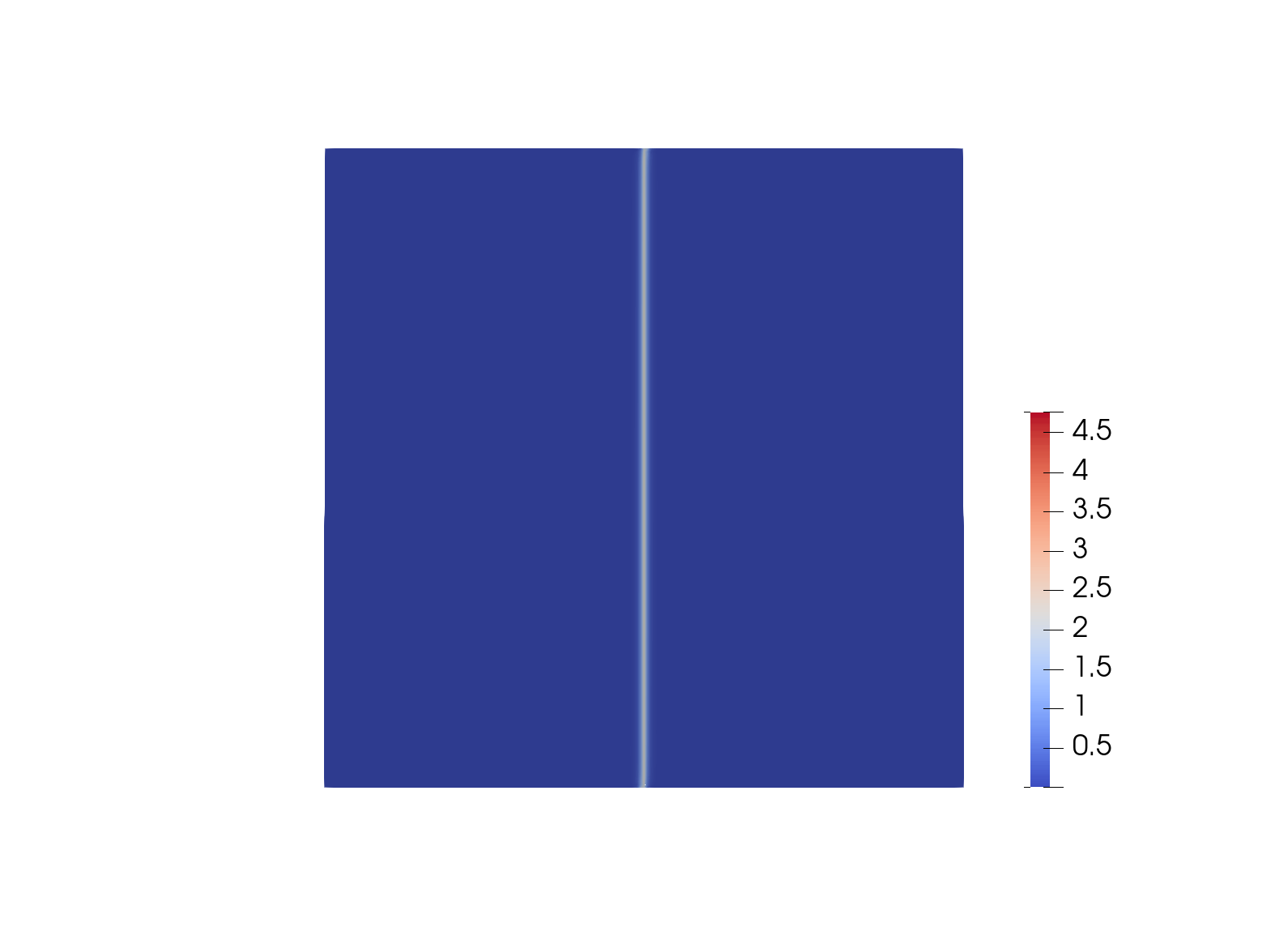}

\caption{Strain fields for both configurations and mesh choices. 
In all cases it is localized on a strip of width $h=0.005.$}\label{fig:strains}
\end{figure}
     
 \end{remark}

It is now instructive to examine a scenario where the solution remains unique but is oriented along a diagonal. This allows us to assess whether mesh independence still holds when the crack direction is oblique and hence more likely to be affected by the mesh geometry. 
To this end, we consider an L-shaped domain with long edge length $L=1$. 

We adopt the same material parameters and boundary conditions as in the previous configuration. Roller supports are applied along the left and bottom edges, while the re-entrant edges at the bottom-left corner are left free to move. Equal perpendicular displacements are prescribed on the right and top edges and are increased uniformly from zero up  to   $U_{\max} = 0.018$ over 6 loading steps.\bl 

\begin{figure}[h]
\centering
\begin{tikzpicture}[scale=4]

\draw[thick] (0.5,0) -- (1,0) -- (1,1) -- (0,1) -- (0,0.5) -- (0.5,0.5) -- cycle;

\foreach \x in {0.1,0.3,0.5,0.7,0.9} {
    \draw[-latex, thick] (\x,1) -- ++(0,0.1);
}

\foreach \y in {0.1,0.3,0.5,0.7,0.9} {
    \draw[-latex, thick] (1,\y) -- ++(0.1,0);
}

\draw[thick] (-0.055,0.5) -- (-0.055,1);     
\draw[thick] (0.5,-0.055) -- (1,-0.055);     

\foreach \x in {0.6,0.75,0.9} {
    \draw[thick] (\x,-0.03) circle (0.02);
}

\foreach \y in {0.6,0.75,0.9} {
    \draw[thick] (-0.03,\y) circle (0.02);
}

\foreach \x in {0.58,0.64,...,0.98} {
    \draw[thick] (\x,-0.1) -- (\x-0.04,-0.06);
}

\foreach \y in {0.58,0.64,...,0.98} {
    \draw[thick] (-0.1,\y) -- (-0.06,\y-0.04);
}

\node[scale=0.9] at (0.5,1.17) {$U_y$};
\node[scale=0.9] at (1.17,0.5) {$U_x$};

\end{tikzpicture}
\caption{Set-up of the numerical simulation.}
\label{fig:L}
\end{figure}
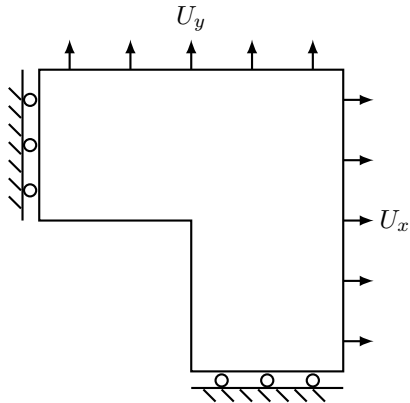

\begin{figure}
\centering
\includegraphics[width=.4\textwidth]{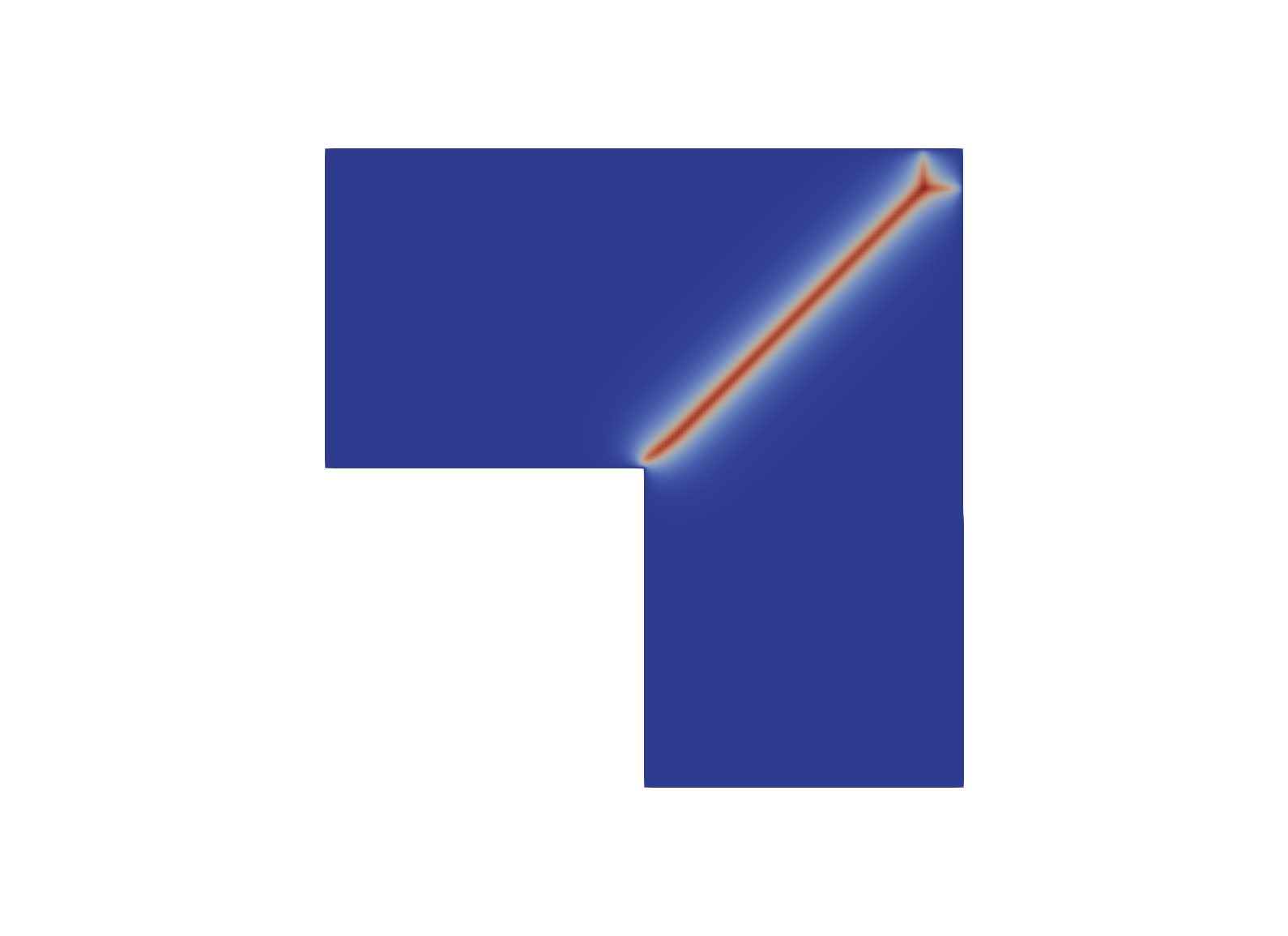}\includegraphics[width=.4\textwidth]{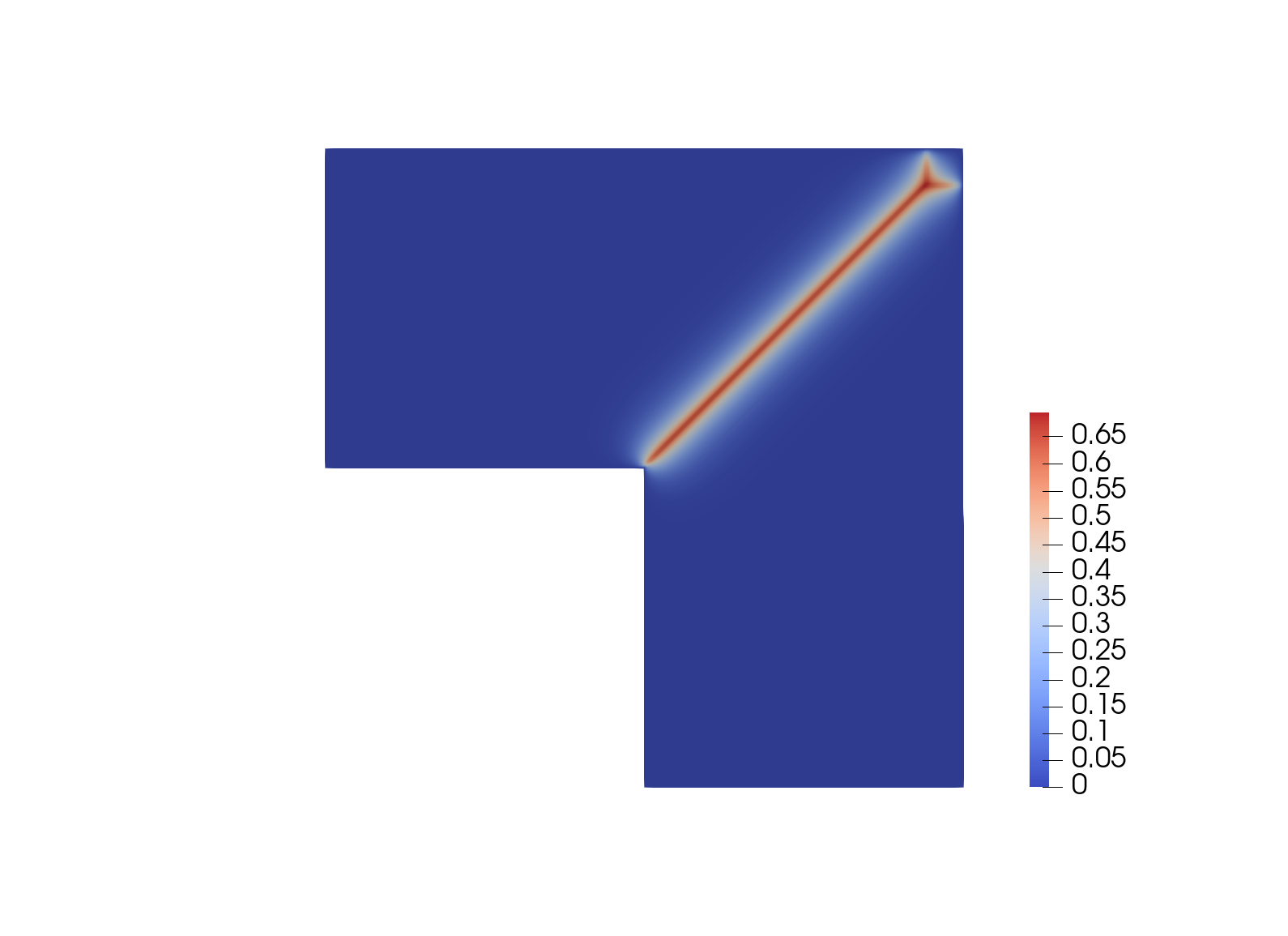}
\caption{Phase fields for mesh $\mathcal T_A$ (left) and mesh $\mathcal T_B$ (left) at $t=39$. In both cases it is localized on a strip of width $\epsilon_h=0.025.$}\label{fig:Lpf}
\end{figure}

\begin{figure}
\centering
\includegraphics[width=.4\textwidth]{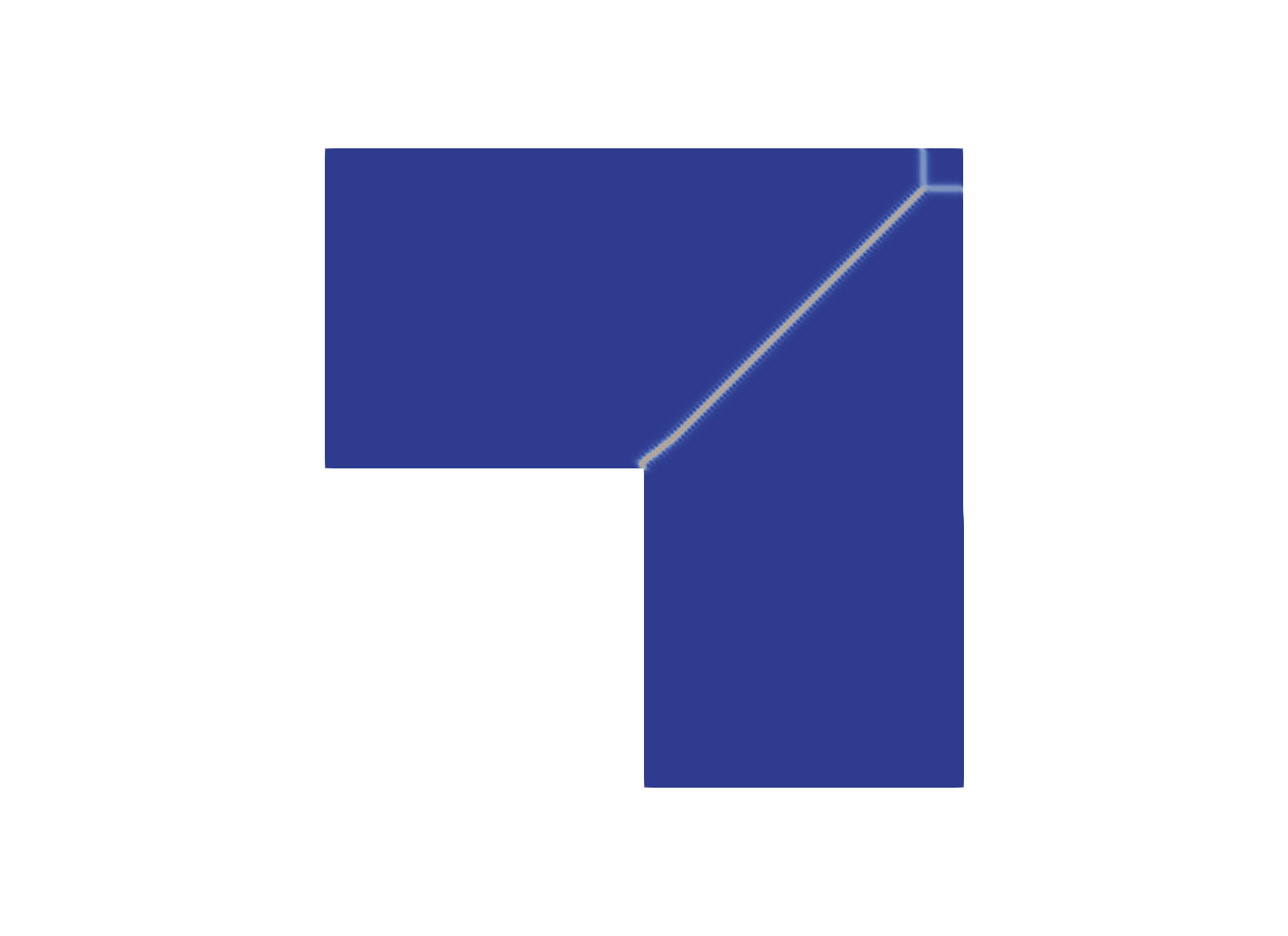}\includegraphics[width=.4\textwidth]{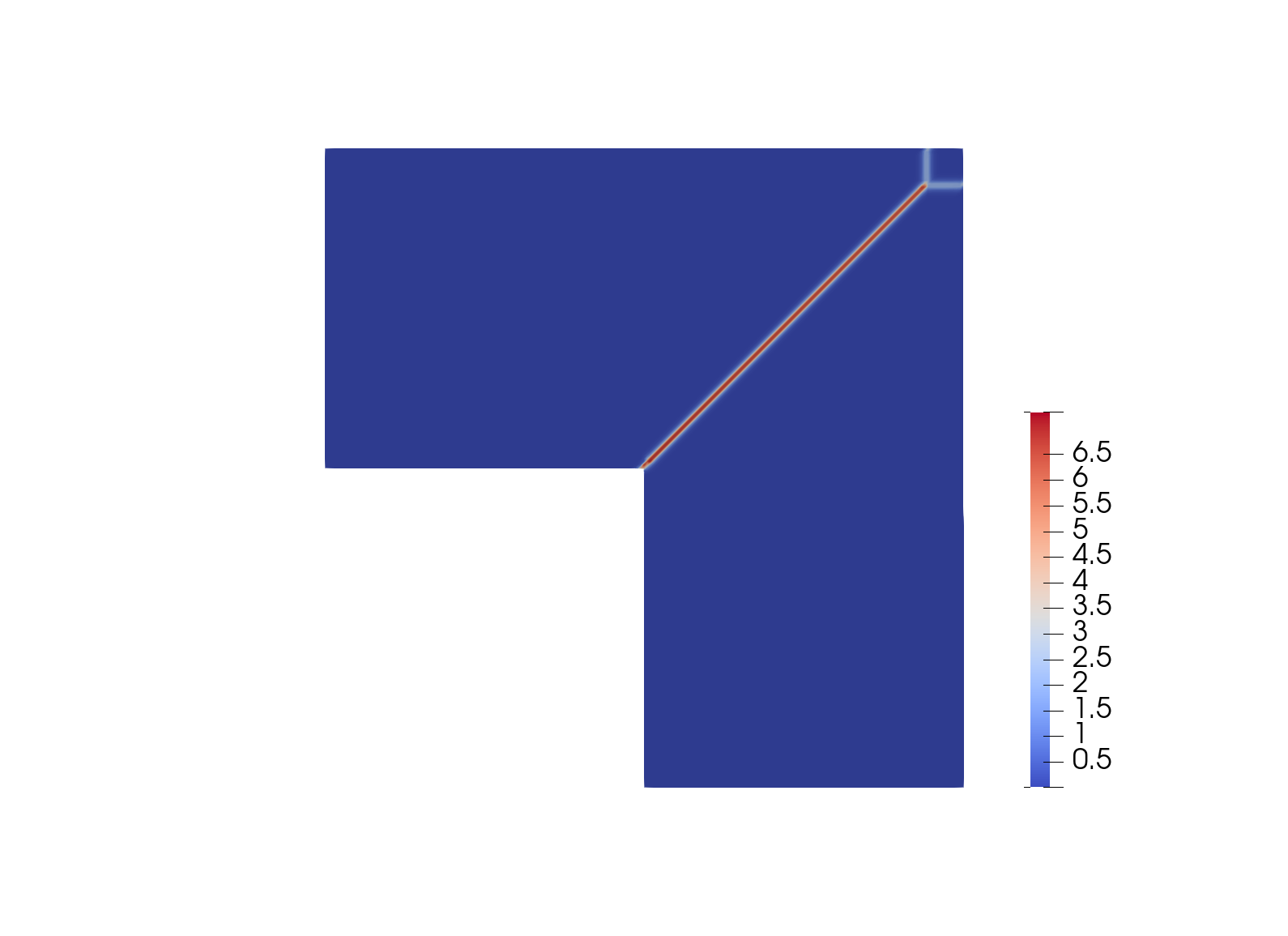}
\caption{Strain for mesh $\mathcal T_A$ (left) and mesh $\mathcal T_B$ (left) at $t=39$. In both cases it is localized on a strip of width $h=0.005.$}\label{fig:strainL}
\end{figure}

\begin{figure} 
\centering
\includegraphics[width=.45\textwidth]{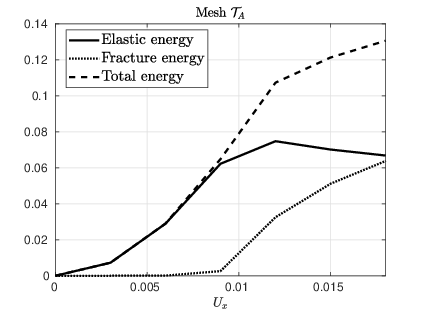}
\includegraphics[width=.45\textwidth]{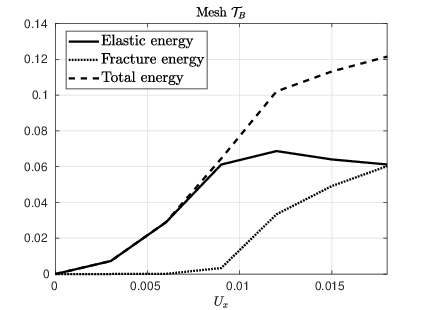}
\includegraphics[width=.45\textwidth]{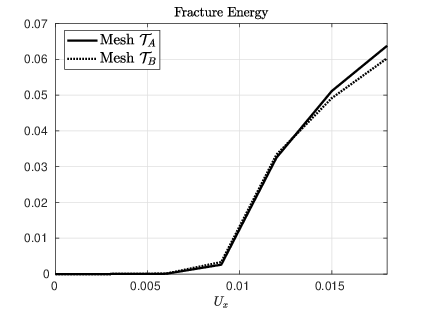}\caption{Evolution of the energies as a function of the displacement in direction $x$. Top row: results obtained using mesh 
$\mathcal T_A$
  (left) and mesh 
$\mathcal T_B$. Bottom row: comparison of fracture energies for the two mesh configurations. 
}\label{fig:cfrEn}\end{figure}

As can be seen in Figure \ref{fig:Lpf}, the crack path is independent of the mesh geometry.
Figure \ref{fig:strainL} depicts the strain fields, that localize in a strip proportional to the mesh width. 
Furthermore, the crack initiates at the same loading step ($U_x=U_y=0.006)$ in both cases and the energy evolutions are similar, as shown in Figure \ref{fig:cfrEn}. 

\newpage

%% file: appendix.tex
\section*{Appendix}
\renewcommand{\thelemma}{{\bf Lemma A.\arabic{lemma}}}
\setcounter{teo}{0}

\renewcommand{\thesection}{\Alph{section}}
\setcounter{section}{1}

\renewcommand{\theequation}{A.\arabic{equation}}
\setcounter{equation}{0}

We present a few technical lemmata on density, relaxation and lower semi-continuity, and properties of the density $f$. Moreover, \ref{subsec:phi_num} contains the numerical setup used to obtain the surface energy density plotted in Figure \ref{fig:phi_1D}.

\subsection{\normalfont Density}

\begin{lemma} \label{lem: denisty} 
    Given $u\in \BV(I)$ such that $\llbracket  u\rrbracket>0$ and $D^cu\geq 0$, there exists a sequence $u_k\in U=\{v\in  \SBV(I)\cap W^{2,\infty}(I\setminus J_{v})\,:\,\#J_{v}<+\infty,\,\llbracket  v\rrbracket>0\}$ such that $u_k\to u$ in $L^1(I)$ and $\limsup_{k\to+\infty}\F(u_k)\leq \F(u)$.
\end{lemma}
\proof 
We construct the sequence 
$u_k$ step by step, first ensuring that it possesses a finite set of jump points.
For each $n\in \mathbb{N}$, let $J_u^n$ denote the set of jump points of $\tilde u$, defined in \eqref{utilde}, with amplitude greater than $\frac1n$:
$$J_u^n=\bigg\{x\in I\,:\,\llbracket  \tilde u(x)\rrbracket>\frac1n\bigg\}.$$
This set has finite cardinality because, since $\phi$ is monotone increasing and $\phi(s) >0$ for $s>0$, the following holds:
$$\F(u)\geq\sum_{J^n_u}\phi(\llbracket  \tilde u\rrbracket)\geq\#J^n_u \phi\bigg(\frac1n\bigg).$$
We call $u'$ the absolutely continuous part of $Du$ and define 
 $$u_n(x)=u(-L)+\int_{[-L,x]}u'\,dx+D^cu([-L,x])+\sum_{J^n_u\cap[-L,x]}\llbracket  \tilde u\rrbracket. $$
Observe that $\llbracket  \tilde u_n\rrbracket>0$, since $\llbracket  \tilde u\rrbracket>0$, $\tilde{u}_n' = u'$, and $D^c \tilde{u}_n = D^c u$. Hence, for every $n\in\mathbb{N}$,  $\F(u_n)\leq \F(u)$  and $$|\F(u_n)-\F(u)|=\sum_{J_{\tilde u}\setminus J_u^n}\phi(\llbracket  \tilde u\rrbracket)\to 0.$$
 Note that 
 $$u(x)-u_n(x)=\sum_{(J_{\tilde u}\setminus J_u^n)\cap[-L,x]}\llbracket  \tilde u\rrbracket\leq \sum_{J_{\tilde u}\setminus J_u^n}\llbracket  \tilde u\rrbracket,$$
so, for every $\epsilon>0$, since $u\in \BV(I)$,
  $$\|u_n-u\|_{L^1(I)}\leq|I|\sum_{J_{\tilde u}\setminus J_u^n}\llbracket  \tilde u\rrbracket<\epsilon$$
 for $n$ sufficiently large. From now on, we fix $\epsilon>0$ and such $n$.
 
 Now, since we want to build a sequence of functions $u_k$ that belong to $\SBV(I)$, we want to get rid of the Cantor part of $u_n$. With this purpose in mind, let us consider a connected component $(a,b)$ of $I\setminus J_u^n$ and uniformly subdivide it, setting:
 $$x_{i,k}=a+i\frac{b-a}{k+1}$$ for $i=0,...,k+1.$ Let us define locally $w(x)=
 D^cu([a,x])$ and 
 $$w^k(x)=\begin{cases}0\qquad \qquad\qquad\,\,\,\quad x\in [a, x_{1,k}]\\
     w(x_{i+1,k})\qquad \qquad x\in (x_{i,k},x_{i+1,k}] \quad \text{for $i=1,...,k$.}
 \end{cases}$$ 
 Note that $w^k (s) = w(s)=0$ and $w^k(b) = w(b)=D^c u ([a,b])$.  
 Of course $\|w-w^k\|_{L^1(a,b)}\to 0$ for $k \to +\infty$, because, being $w$ continuous, it is Cauchy integrable. Moreover, since $\phi(s)\leq \sigma_cs$,
  \begin{equation}\begin{aligned}
      \sum_{J_{w^k}}\phi(\llbracket w^k\rrbracket
  )dx&\leq \sigma_c\bigg(w(x_{2,k})+\sum_{i=2}^{k} w(x_{i+1,k})-w(x_{i,k})\bigg)
  \\&=
\sigma_c\bigg(D^cu([a,x_{2,k}])+\sum_{i=2}^{k} D^cu([x_{i,k},x_{i+1,k}])\bigg)=\sigma_cD^cu([a,b]).
  \end{aligned}\end{equation}
  We now define locally
  
  $$u_{k}(x)=u(a^+)+\int_{[a,x]}u'dx+w^k(x)$$
  and observe that $u_{k}(a)=u(a^+)$ and $u_k(b)=u(b^-)$, hence $\llbracket u_k\rrbracket=\llbracket u_n\rrbracket$ on  $J^n_u$. In particular, the positivity of the jump amplitudes is preserved.
  The set of jump points of $u_k$ is given by $$J_k  =  J_u^n\sqcup J_{w^k} $$ that has cardinality lower than $\#J^n_u+k\cdot(\#J^n_u+1) $.
  We analyze separately the terms of $\F(u_k)$ and we start by observing that: $$\sum_{J_k} \phi(\llbracket u_{k}\rrbracket) \le \sum_{J_u^n}\phi(\llbracket u_{k}\rrbracket) +\sum_{J_{w^k}}\phi(\llbracket u_{k}\rrbracket)\leq \sum_{J_{u}^n}\phi(\llbracket u_{n}\rrbracket)+\sigma_c D^cu(I).$$
  This, together with the fact that $u'_k=u'_n=u'$ and hence $f(u'_k)=f(u')$, leads to:
  $$\F(u_{k})\leq \F(u_n)\leq \F(u).$$
  Finally, for every $\epsilon'>0$
  $$\|u_k-u\|_{L^1(I)}\leq \|u_k-u_n\|_{L^1(I)}+\|u_n-u\|_{L^1(I)}<\|w-w^k\|_{L^1(I)}+\epsilon<\epsilon'$$
  for $k$ sufficiently large.
 A priori, the functions $u_k$ do not necessarily belong to $W^{2,\infty}(I\setminus J_k)$, but we can consider each connected component $(a',b')$ of $I\setminus J_k$ and build locally an approximating sequence of regular functions. In particular we set $u_{k,\lambda}\in W^{2,\infty}(a',b')$ that converges in $H^1(a',b')$ to $u_k|_{(a',b')}$ and such that $u_{k,\lambda}(a')=u_k(a')$ and $u_{k,\lambda}(b')=u_k(b')$. The convergence in $H^1$ guarantees that $\F(u_{k,\lambda})\to \F(u_{k})\leq \F(u)$. By abuse of notation we set $u_k=u_{k,\lambda}$ for $\lambda$ sufficiently large.
\qed

\begin{remark} \label{rem-dens} 
In the proof of \Cref{lem: denisty}, we observed that replacing a Cantor part with a jump discontinuity leads to a lower energy. Consequently, from the perspective of energy minimization, it is convenient to concentrate the singular part of the derivative on jump sets rather than on Cantor-type sets. By a standard truncation argument, in \Cref{lem: denisty} we can assume $\| u_k \|_\infty \le \| g \|_\infty$ in the case $\| u \| \le \|g \|_\infty$. 
\end{remark}

\subsection{\normalfont Lower semi-continuity of the limit} \hfill

\vspace{10pt}
The proof of the lower semi-continuity of $\widetilde\F$ relies on the results presented in \cite{BBBF}. In order to apply them, we introduce the following approximations of the functions $f$ and $\phi$, defined in \eqref{eq:f} and \eqref{eq:phi}:
\begin{equation}
 \label{eq:fn}
f_n(s)=\begin{cases}
    f(s,0) & s\geq  -\frac n{E_0} ,
    \\ 
    \frac{E_0}2n^2-n(s-n) & s \leq -\frac n{E_0} ,
\end{cases}   
\end{equation}
\begin{equation}
    \label{eq:phin}
\phi_n(s)=\begin{cases}
    \phi(s) & s\geq  0 ,
    \\-ns & s<0.
\end{cases}
\end{equation}
These approximations allow us to define the auxiliary functional $\G_n : L^1(I) \to [0, +\infty]$ as follows:
 $$\G_n(u)=\begin{cases}
     { \displaystyle \int_I f_n(u')dx+\sum_{J_{\tilde u}}\phi_n(\llbracket  \tilde u\rrbracket)} & u\in\SBV(I), \\
     +\infty &  u\in L^1(I)\setminus\SBV(I).
 \end{cases}
    $$ 
Adapting the results of \cite{BBBF}, we show that the supremum over $n\in \mathbb{N}$ of the $L^1(I)$-relaxation of $\G_n$ coincides with the functional $\widetilde\F$. To align with the notation used in \cite{BBBF}, we  call $\phi_n^0$ the recession function of $\phi_n$ in the origin, 
$$\phi_n^0(s)
= \lim_{t\to 0^+}\frac{\phi_n(ts)}{t}=\begin{cases}
    \sigma_c s & s\geq 0, \\-ns & s<0, 
\end{cases}$$ 
and $f^\infty_n$ (resp. $g_n^\infty$) the recession function of $f_n$ (resp $g_n$) at infinity:
\begin{equation}
\label{eq: finf}
    f_n^\infty(s)=\limsup_{t\to+\infty}\frac{f_n(ts)}{t}=\begin{cases}
        \sigma_cs & s\geq0, 
        \\-ns & s<0.
    \end{cases}
\end{equation}
The following Lemma summarizes the result of \cite{BBBF} in our setting.
    \begin{lemma}

    The lower semi-continuous envelope of $\G_n$ has the form:
\begin{equation} 
    \bar \G_n(u)=\begin{cases}
        {\displaystyle \int_I g_n(u')dx+\sum_{J_{\tilde u}}h_n(\llbracket  \tilde u\rrbracket)+\int_I g_n^\infty (dD^cu)} & u\in \BV(I), 
        \\+\infty &  u\in L^1(I)\setminus \BV(I), 
    \end{cases}
\end{equation}
where $g_n$ is the inf-convolution of $f_n$ and $\phi_n^0$:
$$g_n(s)= f_n\triangledown\phi^0_n(s)=\min\{f_n(s-r)+\phi^0_n(r)\,:\,r\in \mathbb{R}\}$$ and \begin{equation}\label{hn}h_n(s)=\min\bigg\{H_n(v):=\int_{(0,1)}f^\infty_n(v')dx+\sum_{J_v}\phi_n(\llbracket  v\rrbracket)\,:\,v\in \SBV(0,1),\,v(0)=0,\,v(1)=s\bigg\}.\end{equation}

    \end{lemma}
    \proof In \cite{BBBF}, Theorem 2.13 is stated for functions \( f \) and \( \varphi \) that depend explicitly on the spatial variable \( x \in I \). Additionally, the function \( \varphi \) depends on the normal \( \nu \) to the jump set. In the one-dimensional setting, the normal vector \( \nu(x) \) takes values in \( \{\pm 1\} \), and the scalar product \( \llbracket \tilde{u} \rrbracket(x) \cdot \nu(x) \) reduces to the jump \( \tilde{u}(x^+) - \tilde{u}(x^-) \). Accordingly, we define
\[
\varphi_n(s, \nu) := \phi_n(s \cdot \nu),
\]
so that the functions \( f_n \) and \( \varphi_n \) defined in \eqref{eq:fn} and \eqref{eq:phin}, correspond, respectively, to \( f \) and \( \varphi \) in the setting of \cite{BBBF}. In our notation, we emphasize the dependence on the parameter \( n \in \mathbb{N} \) and instead omit the dependence on $x\in I$, as the functions are independent of the spatial variable. 
It is straightforward to check that $f_n$ and $\varphi_n$ meet the hypothesis (H0)-(H7) of \cite{BBBF}.
Therefore, by the aforementioned Theorem, we obtain an integral representation of the relaxation of $\G_n$ in $\BV(I)$ with respect to the $\BV(I)$-weak topology, namely:
$$\bar\G_n(u)=\inf\bigg\{\liminf_{k\to +\infty}\G_n(u_k)\,:\,u_k\in \SBV(I),u_k\to u \text{ in }L^1(I),\,\sup_k|Du_k|(I)<+\infty\bigg\}.$$
Actually, as observed in \cite{BBBF}, since $f_n\geq 0$, we obtain the same relaxation of $\G_n$ with respect
to the $L^1(I)$ topology.\qed
\begin{lemma}\label{lem: fn phin}The functions $g_n$ and $h_n$ that appear in the integral representation of $\bar\G_n$ are such that:
\begin{equation}
    \label{eq:g}  g_n(s)
=f_n(s),
\end{equation}
\begin{equation}
    \label{eq:h} h_n(s)=\phi_n(s).
\end{equation}
\end{lemma} 
\proof Let's start by proving \eqref{eq:g}, i.e. that for any $r\in \mathbb{R}$,
$f_n(s-r)+\phi_n^0(r)\geq f_n(s)$. From the convexity of $f_n$, we have 
$$ f_n(s-r)\geq f_n(s)+f_n'(s)(-r) $$
which proves the estimate, since $\phi_n^0(r)\geq f_n'(s)r$. Indeed
$$
\phi_n^0(r)=\begin{cases}
    \sigma_c r & r\geq 0, \\-nr & r<0, 
\end{cases} \qquad \qquad 
f_n'(s)=\begin{cases}\sigma_c & s>\frac{\sigma_c}{E_0}, 
\\E_0 s & s\in(-\frac{n}{E_0}, \frac{\sigma_c}{E_0}), \\
-n & s\leq -\frac{n}{E_0}. 
\end{cases}$$

We now prove equation \eqref{eq:h}, recalling the definition of $H_n$ given in \eqref{hn}.
The monotone envelope $\hat v$ of any function $v\in \SBV(0,1),\,v(0)=0,\,v(1)=s$ is such that $H_n(\hat v)\leq H_n(v)$. Let us assume $s>0$, the case $s<0$ being similar. As a consequence the minimum of the functional $H_n$ is attained for a non-decreasing function that takes values between $0$ and $s$. Let us call this function $v$.
We now prove that $v$ has at most one jump point and to do so, let us assume by contradiction that $N=\#J_v>1$. 
Since $\phi_n$ is strictly concave and hence sub-additive it follows that:
$$
\sum_{i=1}^N\phi_n(\llbracket  v\rrbracket)> \phi_n\bigg(\sum_{i=1}^N\llbracket  v\rrbracket\bigg). 
$$
We introduce a function $\bar v\in \SBV(0,1)$ such that $\bar v(0)=0$, $\bar v(1)=s$ and $D\bar v=v'dx+\big(\sum_{i=1}^N\llbracket  v\rrbracket\big)\delta_{x_0}$, where $v'$ is the absolutely continuous part of $Dv$ and $x_0\in(0,1)$. Then
$$H_n(v)=\int_{(0,1)}f^\infty_n(v')dx+\sum_{J_v}\phi_n(\llbracket  v\rrbracket)> 
\int_{(0,1)}f^\infty_n(v')dx+\phi_n\bigg(\sum_{i=1}^N\llbracket  v\rrbracket\bigg)=
H_n(\bar v),$$
 contradicting the minimality of $v$. 
We now prove that the function that minimizes $H_n$ is piecewise constant and has one jump point of amplitude $s$.
Set $0\leq a <b\leq 1$, we assume by contradiction that $v|_{(a,b)}$ is a strictly increasing continuous function. 
This function minimizes $H_n|_{(a,b)}$ 
among the $\SBV$ functions such that $w(a)=v^+(a),\,w(b)=v^-(b)$.
By definition \eqref{eq: finf}, 
$f^\infty_n(v') = \sigma_cv'$ and hence:
\begin{align}
 H_n|_{(a,b)}(v)=\int_a^bf^\infty_n(v')dx=\sigma_c(v^-(b)-v^+(a)).
\end{align}
Since $\sigma_c (v^-(b)-v^+(a))>\phi_n\big( v^-(b)-v^+(a)\big) $, then for a given $x_0\in (a,b)$, $$\tilde v(x)=v^+(a)+(v^-(b)-v^+(a))\chi_{(x_0,b)}(x)$$  is such that $H_n|_{(a,b)}(\tilde v)=\phi_n(v^-(b)-v^+(a))< H_n|_{(a,b)}(v)$, contradicting the minimality of $v|_{(a,b)}$. 
 Hence, substituting this function into \eqref{hn}, we obtain that $h_n(s)=\phi_n(s)$ and the proof is concluded.\qed
\begin{cor}\label{relaxGn}
   For every $n\in \mathbb{N}$, the relaxation of $\G_n$ in $L^1(I)$ has the following form:
    \begin{equation}
        \label{eq: Gn}
    \bar \G_n(u)=\begin{cases}
        {\displaystyle \int_If_n(u')dx+\sum_{J_{\tilde u}}\phi_n(\llbracket  \tilde u\rrbracket)+\int_I f_n^\infty (dD^cu)} & \text{if $u\in \BV(I)$,} \\+\infty &  \text{otherwise.}      
    \end{cases}
    \end{equation} 
\end{cor}

\begin{cor}\label{supGn} It holds $\widetilde \F=\sup_n \bar \G_n$, in particular $\widetilde{\mathcal{F}}$ is lower semi-continuous.
\end{cor}

\proof
 We substitute definitions \eqref{eq:fn}, \eqref{eq:phin} and \eqref{eq: finf} into \eqref{eq: Gn}; note that the integrals appearing in \eqref{eq: Gn} are defined on mutually disjoint subsets of $\Omega$, moreover the densities $f_n$, $\phi_n$ and $f^\infty_n$ are monotone increasing, with respect to $n$, and such that
 $$
      \sup_n f_n (s,0) = f (s,0), \qquad \sup_n \phi_n(s) = \begin{cases} \phi(s) & s \ge 0 \\ + \infty & s < 0 . \end{cases} 
      \qquad 
      \sup_n f^\infty_n (s) = 
      \begin{cases}
          \sigma_c s & s \ge 0 \\ + \infty & s < 0 . 
      \end{cases} 
$$ Then, we get 
$$\sup_n \bar \G_n(u)=\begin{cases} {\displaystyle \int_If(u',0)dx+\sum_{J_{\tilde u}}\phi(\llbracket  \tilde u\rrbracket)+\int_I \sigma_c (dD^cu) } & \text{if $u\in \BV(I),\,\llbracket \tilde u\rrbracket> 0$,\, $D^cu\geq 0$}\\+\infty &  \text{otherwise,}      
\end{cases}$$ that is the definition \eqref{eq: F} of $ \widetilde\F(u)$. 

\begin{remark}
   Given $u\in \BV(I)$ such that $\llbracket  u\rrbracket>0$ and $D^cu\geq 0$, let $u_k\in U$ 
     be the sequence defined in \Cref{lem: denisty}. By lower semi-continuity of $\F$ it follows that $\lim_{k\to+\infty}\F(u_k)=  \F(u)$.
\end{remark}



\smallskip

\subsection{\normalfont Properties of energy density $f$} \hfill

\medskip
\noindent We provide first this concentration lemma. 

\begin{lemma}\label{fBar}
Let $f$ be defined in \eqref{eq:f}, and introduce the threshold value $s(r)=a(r)\frac{\sigma_c}{E_0}$.
  Let $r_i\in[0,1]$ for $i=1,...,n$, where $r_n=\max_{i=1,..,n}r_i$. Let $s_i \ge s(r_n)$ for $i=1,...,n$, then the following inequality holds:
  $$\sum_{i=1}^nf(s_i,r_i)\geq \sum_{i=1}^{n-1}f(s(r_n),r_i)+f\bigg(\sum_{i=1}^{n}s_i-(n-1)s(r_n),r_n\bigg).$$
\end{lemma}

\proof
Set $S=( \sum_{i=1}^{n}s_i)-(n-1)s(r_n)=s_n+\sum_{i=1}^{n-1}(s_i-s(r_n))
$, 
the thesis is equivalent to proving that:
 $$\sum_{i=1}^{n-1}\bigg(f(s_i,r_i)-f(s(r_n),r_i)\bigg)\geq f(S,r_n)-f(s_n,r_n).$$
Since $S \ge s_n\geq s(r_n)$ and hence $f|_{[s_n,S]}(\cdot,r_n)$ is linear, the right hand side is:
$$f(S,r_n)-f(s_n,r_n)=\int_{s_n}^{S}\partial_sf(s,r_n)ds=a(s(r_n))\sigma_c (S-s_n).$$
On the other hand, since $\partial_sf(s,\cdot)$ is decreasing and $r_i\leq r_n$:
 $$f(s_i,r_i)-f(s(r_n),r_i)=\int_{s(r_n)}^{s_i}\partial_sf(s,r_i)ds \geq \int_{s(r_n)}^{s_i}\partial_sf(s,r_n)ds=a(s(r_n))\sigma_c(s_i-s(r_n)).$$
 It follows that
 $$\sum_{i=1}^{n-1}\bigg(f(s_i,r_i)-f(s(r_n),r_i)\bigg)\geq a(s(r_n))\sigma_c \sum_{i=1}^{n-1} (s_i-s(r_n))=a(s(r_n))\sigma_c(S-s_n)$$
 and the Lemma is proven. \qed

\begin{lemma}\label{l.fLip} Let $f$ be defined in \eqref{eq:f}, then for $r \in [0,1]$ we have 
$$
    f(s,0) - f(s,r) \le C | r| | s| .
$$
\end{lemma}

\proof Let $s(r) =  a(r)  \frac{\sigma_c}{E_0} $ be the threshold appearing in the definition of $f$. Clearly $0 \le s(r) \le s(0)$. For $s \le s(r)$ we have $f(s,0) = f(s,r)$ and there is nothing to prove. For $s > s(r)$ let us write 
$$
    f(s,0) - f(s,r) 
    = \int_{s(r)}^{s} f' (z , 0) - f' ( s(r) , r) \, dz =  \int_{s(r)}^{s} f' (z , 0) - f' ( s(r) , 0) \, dz.
$$
Hence, for $s (r) \le s \le s(0)$ we have
\begin{align}
    f(s,0) - f(s,r) 
    & 
    = \int_{s(r)}^{s} E_0 ( z - s(r) ) \, dz 
    = \tfrac12 E_0 ( s - s(r) )^2 \le \tfrac12 E_0 ( s - s(r) ) ( s(0) - s(r) ) \\ & = \tfrac12 \sigma_c ( s - s(r)) ( a(0) - a(r) ) \le C (s - s(r)) |r| \le C |s| |r| . 
\end{align}
For $s > s(0)$ let us write 
$$
    f( s, 0) - f (s, r) = \int_{s(r)}^{s(0)} f' (z , 0) - f' ( s(r) , 0) \, dz + \int_{s(0)}^{s} f' (z , 0) - f' ( s(r) , 0) \, dz .
$$
The first integral is estimated by $C (s (0) - s (r)) |r| $ (see above). For the second it is enough to write
\begin{align}
    \int_{s(0)}^{s} f' (z , 0) - f' ( s(r) , 0) \, dz 
    & = \int_{s(0)}^{s} E_0 ( s(0) - s(r)) \, dz = E_0 ( s(0) - s(r) ) ( s - s(0)) \\
    & \le \sigma_c ( a (0) - a(r) ) ( s - s(0)) \le C |r| ( s - s (0)) .
\end{align}
Joining the inequalities gives 
$$
    f( s, 0) - f (s, r) \le C |r| ( s - s(r) ) \le C |r| |s| , 
$$
which concludes the proof.  \qed

\subsection{Numerical surface energy density for the 1D model}
 \label{subsec:phi_num}

\begin{figure}[H]
    \centering
    \includegraphics[scale=1]{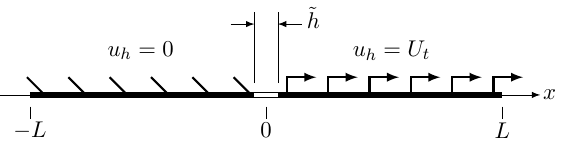}
    \caption{Setup for the 1D surface energy density test. The white region in the bar scheme represents the tiny element (size $\tilde{h}=h/25$, not to scale) dividing the bar into three regions.}
    \label{fig:phi_1D_setup}
\end{figure}

We compute the numerical surface energy density in Fig.~\ref{fig:phi_1D} using the one-dimensional \texttt{FEniCSX} finite element implementation described in Section 5.2 of \cite{vicentini2024energy}, applied to the problem illustrated in Fig.~\ref{fig:phi_1D_setup}. We set $L = 1$, $E_0 = 10^4$, $G_c = 10^{-3}$, $\sigma_c = 5$ and $\epsilon_h=0.4$. The element size $h$ satisfies $\epsilon_h/h \approx 5$ throughout the bar, except for a tiny central element of size $\tilde{h}=h/25$, introduced to more accurately capture the displacement jump. As shown in Fig.~\ref{fig:phi_1D_setup}, this tiny central element divides the bar into three regions. The left region is fixed with $u_h = 0$, whereas the right region undergoes a rigid displacement $U_t$, uniformly increased from $0$ to the maximum value $5\times 10^{-3}$ in 50 time steps. To obtain the numerical dependence of the surface energy density $\phi$ on the jump $j$ in Fig.~\ref{fig:phi_1D}, we consider $j = U_t$.

